\documentclass[12pt]{article}

\usepackage{amsmath,amstext,amsfonts,amsthm,amssymb,eucal,bm,diagrams}

\newcommand{\dL}{\mathbb{L}}

\newcommand{\Hqe}{\operatorname{Hqe}}
\newcommand{\maf}{\operatorname{mf}}

\newcommand{\abs}{\operatorname{abs}}

\newcommand{\be}{\begin{equation}}
\newcommand{\ee}{\end{equation}}

\newcommand{\Per}{\operatorname{Per}}

\newcommand{\bP}{{\mathbf P}}

\newcommand{\Id}{\operatorname{Id}}

\newcommand{\ch}{\operatorname{ch}}

\newcommand{\Com}{\operatorname{Com}}

\newcommand{\Tr}{\operatorname{Tr}}
\newcommand{\gr}{\operatorname{gr}}

\newcommand{\und}{\underline}

\newcommand{\OO}{{\cal O}}

\newcommand{\DD}{{\cal D}}

\newcommand{\tot}{\operatorname{tot}}

\newcommand{\tors}{\operatorname{tors}}

\newcommand{\G}{{\mathbb G}}

\newcommand{\hra}{\hookrightarrow}

\newcommand{\lan}{\langle}
\newcommand{\ran}{\rangle}

\newcommand{\CC}{{\cal C}}

\newcommand{\Proj}{\operatorname{Proj}}
\renewcommand{\P}{{\mathbb P}}

\newcommand{\Pic}{\operatorname{Pic}}

\newcommand{\ga}{\gamma}
\newcommand{\de}{\delta}
\newcommand{\eps}{\epsilon}

\newcommand{\im}{\operatorname{im}}

\numberwithin{equation}{section}

\newtheorem{theor}{Theorem}[subsection]
\newtheorem{thm}[theor]{Theorem}
\newtheorem{lem}[theor]{Lemma}

\newtheorem{prop}[theor]{Proposition}

\newtheorem{cor}[theor]{Corollary}
{  \theoremstyle{definition}
           \newtheorem{defi}[theor]{Definition}
           \newtheorem{rem}[theor]{Remark}

}

\newcommand{\Pf}{\noindent {\it Proof}}
\newcommand{\id}{\operatorname{id}}

\newcommand{\ov}{\overline}
\newcommand{\we}{\wedge}

\newcommand{\rk}{\operatorname{rk}}

\renewcommand{\AA}{{\cal A}}

\newcommand{\FF}{{\cal F}}

\newcommand{\TT}{{\cal T}}
\newcommand{\XX}{{\cal X}}
\newcommand{\HH}{{\cal H}}
\newcommand{\PP}{{\cal P}}

\renewcommand{\SS}{{\cal S}}

\newcommand{\LL}{{\cal L}}

\newcommand{\st}{\operatorname{st}}
\newcommand{\Du}{{\mathbb D}}
\newcommand{\Om}{\Omega}

\newcommand{\Hom}{\operatorname{Hom}}

\newcommand{\Ext}{\operatorname{Ext}}

\renewcommand{\a}{\alpha}
\renewcommand{\b}{\beta}

\newcommand{\om}{\omega}
\newcommand{\De}{\Delta}
\newcommand{\la}{\lambda}
\renewcommand{\th}{\theta}
\newcommand{\C}{{\mathbb C}}

\newcommand{\A}{{\mathbb A}}

\newcommand{\Z}{{\mathbb Z}}
\newcommand{\Q}{{\mathbb Q}}
\newcommand{\La}{\Lambda}
\newcommand{\Ga}{\Gamma}
\newcommand{\Td}{\operatorname{Td}}

\newcommand{\wt}{\widetilde}
\newcommand{\ot}{\otimes}

\newcommand{\sub}{\subset}
\newcommand{\com}{\operatorname{com}}
\newcommand{\ed}{\qed\vspace{3mm}}

\newcommand{\Qcoh}{\operatorname{Qcoh}}

\newcommand{\tr}{\operatorname{tr}}
\newcommand{\qgr}{\operatorname{qgr}}

\newcommand{\pa}{\partial}

\newcommand{\Sg}{\operatorname{Sg}}

\newcommand{\At}{\operatorname{At}}

\newcommand{\MF}{{\operatorname{MF}}}

\newcommand{\DMF}{\operatorname{DMF}}

\providecommand{\keywords}[1]{\textit{Keywords.} #1}
\providecommand{\MSC}[1]{2010 \textit{Mathematics Subject Classification.} #1}

\title{Homogeneity of cohomology classes associated with Koszul matrix factorizations}
\author{Alexander Polishchuk}

\date{}

\begin{document}
\maketitle
\begin{abstract} In this work we prove the so called {\it dimension property} for the cohomological field theory
associated with a homogeneous polynomial $W$ with an isolated singularity, in the algebraic framework of \cite{PV-CohFT}.
This amounts to showing that some cohomology classes on the Deligne-Mumford moduli spaces of stable curves,
constructed using Fourier-Mukai type functors associated with matrix factorizations, live in prescribed dimension.
The proof is based on a homogeneity result established in \cite{PV-Wc} for certain characteristic classes of
Koszul matrix factorizations of $0$. To reduce to this result we use the theory of Fourier-Mukai type functors involving
matrix factorizations and the natural rational lattices in the relevant Hochschild homology spaces, as well as 
a version of Hodge-Riemann bilinear relations for Hochschild homology of matrix factorizations.
Our approach also gives a proof of the dimension property for the cohomological field theories associated with some
quasihomogeneous polynomials with an isolated singularity.
\end{abstract}

\MSC{Primary  81T70; Secondary 14B05, 16E40}


\section{Introduction}

\subsection{Dimension property in Fan-Jarvis-Ruan-Witten theory}\label{dim-sec-1}

Let $W\in\C[x_1,\ldots,x_n]$ be a quasihomogeneous polynomial with an isolated singularity at the origin.
Fan, Jarvis and Ruan introduced in \cite{FJR} an analog of the Gromov-Witten theory associated with $W$ and with
a finite subgroup $G\sub (\C^*)^n$ of diagonal symmetries of $W$ (such that $G$ contains the exponential grading operator
$J$ associated with the weights of the variables $x_1,\ldots,x_n$).
This theory, often referred to as Fan-Jarvis-Ruan-Witten theory ({\it FJRW-theory}), consists of a collection of maps
$$\La_{g,r}:H^{\ot r}\to H^*(\ov{M}_{g,r},\C),$$
where $\ov{M}_{g,r}$ is the Deligne-Mumford compactification of the moduli spaces of curves with $r$ marked
points, and $H=H_{W,G}$ is a finite-dimensional vector space associated with $(W,G)$ (called the {\it space state} of
the theory).
The maps $\La_{g,r}$ satisfy some gluing axioms on the boundary components of $\ov{M}_{g,r}$, 
that constitute the notion of a {\it cohomological field theory}, introduced by Kontsevich and Manin
\cite{KM}. In fact, the state space of the FJRW-theory has a decomposition
$$H_{W,G}=\bigoplus_{\ga\in G}H(W_\ga)^G,$$
where $W_\ga=W|_{(\A^n)^\ga}$, the restriction of $W$ to the space of $\ga$-invariants, and 
\begin{equation}\label{H-W-eq}
H(W(x_1,\ldots,x_n)):=(\Om_{\A^n}^n/(dW\we \Om_{\A^n}^{n-1}))
\end{equation}
(the latter definition is applied to all $W_\ga$).
Each component of the map $\La_{g,r}$ factors as a composition
$$H(W_{\ga_1})^G\ot\ldots\ot H(W_{\ga_r})^G\rTo{\phi_g(\ga_1,\ldots,\ga_r)} 
H^*(\SS_{g,G}(\ga_1,\ldots,\ga_r),\C)\to H^*(\ov{M}_{g,r},\C),$$
where $\SS_{g,G}(\ga_1,\ldots,\ga_r)\to \ov{M}_{g,r}$ is some finite covering, corresponding to choices of generalized
spin-structure (of type $\ga_1,\ldots,\ga_r$) on a curve. For details, see \cite{FJR} and \cite{PV-CohFT}.

Whereas in the original approach of \cite{FJR} the maps $\phi_g(\ga_1,\ldots,\ga_r)$ were defined by studying a certain
PDE (Witten's equation), in \cite{PV-CohFT} we constructed these maps using Hochschild homology and the
categories of matrix factorizations. More precisely, we use natural embeddings 
$$H(W)^G\sub HH_*(\MF_G(W)),$$
where $\MF_G(W)$ is the category of $G$-equivariant matrix factorizations of $W$, and construct the maps 
$\phi_g(\ga_1,\ldots,\ga_r)$ as maps induced on Hochschild homology by some Fourier-Mukai type functor
$$\MF_G(W_{\ga_1}\oplus\ldots\oplus W_{\ga_r})\to D^b(\SS_{g,G}(\ga_1,\ldots,\ga_r)),$$
where $D^b(X)$ denotes the derived category of coherent sheaves on $X$.

Conjecturally, the algebraic approach of \cite{PV-CohFT} produces the same theory as in \cite{FJR}, however, this is currently known to be true only for simple singularities (see \cite[Sec.\ 7]{PV-CohFT}),
in the so-called narrow sectors (see \cite{CLL}), and for most invertible polynomials and the maximal groups of symmetries
(see \cite{Guere}). In general, the hope is that the algebraic approach will be more accessible for calculations (as for example,
the work \cite{Guere} indicates), so it is important to establish algebraically all the properties of the FJRW-theory.

One of the properties of the maps $\phi_g(\ga_1,\ldots,\ga_r)$ which arises naturally in the analytic approach of \cite{FJR} is 
the {\it dimension property} stating that
\begin{equation}\label{dimension-axiom}
\im(\phi_g(\ga_1,\ldots,\ga_r))\sub H^{2D_g(\ga_1,\ldots,\ga_r)+n_1+\ldots+n_r}(\SS_{g,G}(\ga_1,\ldots,\ga_r),\C),
\end{equation}
with $n_i=\dim (\A^N)^{\ga_i}$ and 
$$D_g(\ga_1,\ldots,\ga_r)=(g-1)\hat{c}+\iota_{\ga_1}+\ldots+\iota_{\ga_r}=-\sum_{j=1}^n\chi(C,L_j),$$
where $(C,L_1,\ldots,L_n)$ is a smooth curve with a generalized spin-structure 
from the moduli space $\SS_{g,G}(\ga_1,\ldots,\ga_r)$ and the numbers $\hat{c}=\hat{c}_W$, 
$\iota_{\ga_1},\ldots,\iota_{\ga_r}$ are
determined using the weights of the variables $x_1,\ldots,x_n$ (see \cite[Sec.\ 3.2]{FJR}).

This property is not at all clear in the algebraic framework of \cite{PV-CohFT}.
The goal of this paper is to prove 
the dimension property in this framework
assuming that $W(x_1,\ldots,x_n)$ is a homogeneous polynomial, i.e., the degrees of the variables
are $\deg(x_1)=\ldots=\deg(x_n)=1$. 

More generally, for a quasihomogeneous polynomial $W(x_1,\ldots,x_n)$, where $\deg(x_i)=d_i>0$, 
we can define a homogeneous polynomial
\begin{equation}\label{quasihom-modif-eq}
\wt{W}(y_1,\ldots,y_n)=W(y_1^{d_1},\ldots,y_n^{d_n})
\end{equation}
in new variables $y_1,\ldots,y_n$ with $\deg(y_i)=1$.
We will prove the dimension property for the algebraic cohomological field theory associated with $(W,G)$ (for any $G$)
provided $\wt{W}$ still has an isolated singularity at $0$.

\begin{thm}\label{dim-property-thm}
Let $W(x_1,\ldots,x_n)$ be a quasihomogeneous polynomial with an isolated singularity, $G$ a finite group of diagonal
symmetries of $W$, containing the exponential grading element. Assume that $\wt{W}$ still has an isolated singularity at $0$.
Then the maps $\phi_g(\ga_1,\ldots,\ga_r)$ defined in \cite{PV-CohFT} satisfy the dimension property
\eqref{dimension-axiom}.
\end{thm}

This will be deduced from a more general Theorem \ref{main-thm} formulated below.

\subsection{Purity of dimension for functors associated with Koszul matrix factorizations}\label{purity-sec}

Let $W(x_1,\ldots,x_n)$ be a quasihomogeneous polynomial of degree $d$, where $\deg(x_i)=d_i>0$,
with an isolated singularity at the origin,
and denote by $\MF_{\G_m}(W)$ the category of $\G_m$-equivariant matrix factorizations of $W$ (see Sec.\ \ref{mf-sec}
below). Here $\G_m$ acts on $\A^n$ by 
\begin{equation}\label{weighted-action-eq}
\la\cdot (x_1,\ldots,x_n)=(\la^{d_1}x_1,\ldots,\la^{d_n}x_n).
\end{equation}

Let also $X$ be a smooth projective variety.
We are going to prove a certain purity of dimension for the maps 
$$\HH(W):=HH_*(\MF_{\G_m}(W))\to H^*(X,\C)$$
induced by Fourier-Mukai functors
$$\DMF_{\G_m}(W)\to D^b(X)$$
of a special kind, where $\DMF_{\G_m}(W)$ is the derived category of $\G_m$-equivariant matrix factorizations of $W$.

Here is the precise setup. Assume that $A$ is a $\G_m$-equivariant vector bundle on $X$, where
$\G_m$ acts trivially on $X$, 
equipped with a surjective $\G_m$-morphism of $\OO_X$-modules 
$$z:A\to \OO_X^n,$$
where $\G_m$ acts on $\OO_X^n$ with the weights $(d_1,\ldots,d_n)$
Let $\tot(A)$ be the total space of this vector bundle, and let $p:\tot(A)\to X$ be the natural projection.
Note that $z$ corresponds to a morphism
$$Z:\tot(A)\to \A^n,$$
linear on the fibers of $p$ and $\G_m$-equivariant with respect to the action \eqref{weighted-action-eq} on $\A^n$.
Let $B$ be another $\G_m$-equivariant vector bundle on $X$ and suppose we have $\G_m$-morphisms of $\OO_X$-modules
$$\a: \oplus_{i=1}^N S^i(A)\to B^\vee\{d\}, \ \ \b:A\to B,$$
where $S^m(\cdot)$ denotes the $m$th symmetric power, and $\{d\}$ denotes the twist by the character $\la\mapsto \la^d$ of
$\G_m$
We can view $\a$ and $\b$ as $\G_m$-invariant sections of induced bundles on $\tot(A)$:
$$\a\in H^0(\tot(A),p^*B^\vee\{d\}), \ \ \b\in H^0(\tot(A),p^*B).$$
The main assumption is that these sections satisfy
$$\lan \a,\b\ran=-Z^*W,$$
and that the common vanishing locus of $\a$ and $\b$ coincides with the zero section in $\tot(A)$.
Then we have a $\G_m$-equivariant Koszul matrix factorization $\{\a,\b\}$ of $-Z^*W$ on $\tot(A)$, supported at the zero section 
 (see Sec.\ \ref{mf-sec}).
We can use this matrix factorization and the diagram
\begin{diagram}
&&\tot(A)\\
&\ldTo{Z}&&\rdTo{p}&\\
\A^n&&&& X
\end{diagram}
to define a Fourier-Mukai type functor 
\begin{equation}\label{main-Phi-eq}
\Phi:\DMF_{\G_m}(\A^n,W)\to D^b(X): E\mapsto p_*(Z^*E\ot \{\a,\b\})
\end{equation}
(more precisely, this is the functor $\Phi_P$ for $P=\{\a,\b\}$ defined in Sec.\ \ref{mf-ker-sec} below).

This functor has a natural realization on the dg-level which in particular allows
to consider the induced map on the Hochschild homology 
$$\phi=\Phi_*:\HH(W)\to HH_*(X).$$
The Hochschild-Kostant-Rosenberg isomorphism together with the Hodge theory give an
identification $HH_*(X)\simeq H^*(X,\C)$ (see Sec.\ \ref{Chern-sec}).

\begin{thm}\label{main-thm} 
In the above situation assume in addition that the homogeneous polynomial
$\wt{W}$ given by \eqref{quasihom-modif-eq} still has an isolated singularity. 
Then for any $x\in \HH(W)$ one has
$$\Td(A)^{-1}\Td(B)\phi(x)\in H^D(X,\C)\sub H^*(X,\C),$$
where $\Td(\cdot)$ denotes the Todd class, and
\begin{equation}\label{D-eq}
D=2\rk B-2\rk A+n.
\end{equation}
\end{thm}

\subsection{Outline of the proof}

The proof combines some ideas of noncommutative Hodge theory, the relation between graded matrix factorizations
and derived categories of coherent sheaves on hypersurfaces (Orlov's equivalence), and a purity result from
\cite{PV-Wc}.

We start by rewriting the assertion using the left adjoint map to $\phi$ with respect to the canonical pairings on 
the Hochschild homology. Recall that for a smooth and proper dg category $\CC$ the Hochschild homology is equipped
with a canonical nondegenerate pairing $\lan\cdot,\cdot\ran_\CC$, such that the maps on Hochschild homology
induced by an adjoint pair of functors are adjoint with respect to the canonical pairings 
(see Sections \ref{can-pair-sec} and \ref{Hoch-maps-sec} for details).
Thus, the left adjoint map to $\phi$ is given by
$$\psi=\Psi_*:H^*(X)\simeq HH_*(X)\to \HH(W),$$
where $\Psi:D^b(X)\to \DMF_{\G_m}(W)$ is the left adjoint functor to $\Phi$.
Thus, we can rewrite the condition $\a\cdot\phi(\HH(W))\sub H^D(X)$, where $\a=\Td(A)^{-1}\Td(B)$, as 
$$\psi(^{\perp}(\a^{-1}H^D(X)))=0.$$

Under the Hochschild-Kostant-Rosenberg isomorphism $HH_*(X)\simeq H^*(X)$ the canonical pairing takes the form
$$\lan\cdot,\cdot\ran_{D^b(X)}=\int_X \kappa(a)\cdot b\cdot\Td_X,$$
where $a,b\in H^*(X)$, $\kappa$ is the linear operator on $H^*(X)$, such that $\kappa(c)=(-1)^qc$ for $c\in H^{p,q}$
(see Section \ref{can-pair-sec}). This implies that the left orthogonal to $\a^{-1}H^D(X)\sub H^*(X)$ with respect to the
canonical pairing is
$$^{\perp}(\a^{-1}H^D(X))=\bigoplus_{j\neq 2\dim X-D}\kappa(\a\cdot\Td_X^{-1})H^j(X),$$
so we need to check that for each $y\in H^j(X)$, where $j\neq 2\dim X-D$, one has
\begin{equation}\label{psi-vanishing-eq}
\psi\bigl(\kappa(\a\cdot\Td_X^{-1})\cdot y\bigr)=0.
\end{equation}

Next, we recall that the Hochschild homology of the category of matrix factorizations has a canonical decomposition
(see \cite[Thm.\ 2.6.1]{PV-CohFT})
\begin{equation}\label{Hoch-decomp-eq}
\HH(W)=HH_*(\MF_{\G_m}(W))\simeq \bigoplus_{\ga\in \mu_d} H(W_\ga)^{\mu_d},
\end{equation}
where $W_\ga=W|_{(\A^n)^\ga}$, and $H(W)$ is given by \eqref{H-W-eq}. 
In fact, \eqref{Hoch-decomp-eq} is exactly the 
decomposition of $HH_*(\MF_{\G_m}(W))$ into isotypical components with respect to the natural action of $\Z/d$ on it
(see \cite[Thm.\ 2.6.1(ii)]{PV-CohFT}).
Let $\Pi:\HH(W)\to\HH(W)$ denote the projector onto the summand
$H(W)^{\mu_d}$, corresponding to $\ga=1$.
Using the characterization of this summand as $\Z/d$-invariants in $\HH(W)$
we check that $\phi=\phi\Pi$ (see Lemma \ref{phi-inv-lem}), 
and hence the image of $\psi$ is contained in $H(W)^{\mu_d}\sub\HH(W)$.

Now the idea is that $H(W)^{\mu_d}$ should be thought of as an analog of the primitive middle cohomology.
Recall that if $X\sub\P^N$ is a smooth projective variety of dimension $n$ then the primitive part of the middle
cohomology $P^n(X)\sub H^n(X)$ is defined as the kernel of the operator of multiplication with $c_1(H)$, where
$H$ is the hyperplane class on $X$. The classical Hodge-Riemann relations imply (see Lemma \ref{HR-lem}) that for
a nonzero class $a\in P^n_H(X)$ one has 
$$(-1)^{\frac{n(n+1)}{2}}\lan a,\ov{a}\ran_{D^b(X)}>0,$$
where $a\mapsto\ov{a}$ is the complex conjugation associated with the real structure on $H^*(X)$.

The key step of the proof is establishing an analog of this property for the canonical pairing on $H(W)^{\mu_d}\sub\HH(W)$.
One missing piece of structure that we need for this is a real structure on $\HH(W)$. In fact, in general one expects to
have a natural rational lattice in the Hochschild homology of any smooth proper dg-category (see \cite{Kon}, \cite{Blanc}).
In the case of an admissible subcategory $\CC$ of $D^b(Y)$ such a rational lattice can be constructed 
easily using the realization of its Hochschild homology as an image of a rational projector on $H^*(Y)$ (see 
Section \ref{rational-lattice-adm-sub-sec}). Furthermore, the obtained rational lattices are
compatible with the maps on Hochschild homology induced by functors of Fourier-Mukai type.

Using the natural embedding $H(W)^{\mu_d}\sub H(\wt{W})^{\mu_d}$,
where $\wt{W}$ is given by \eqref{quasihom-modif-eq} we reduce the situation to the homogeneous case (where $\deg(x_i)=1$).
For homogeneous $W$, 
we apply Orlov's result, connecting the category of matrix factorizations $\DMF_{\G_m}(W)$ with the derived category 
of the corresponding projective hypersurface $X=(W=0)$, to realize $\DMF_{\G_m}(W)$ as such an admissible subcategory
(see Lemma \ref{MF-admissible-lem}). Then using explicit descriptions of the canonical pairing $\lan\cdot,\cdot\ran_W$
on $\HH(W)$ and of the Chern characters of matrix factorizations in \cite{PV-HRR}, we prove the following
property of the subspace $H(W)^{\mu_d}\sub\HH(W)$, which makes it an analogue of the primitive cohomology:
$H(W)^{\mu_d}$ is orthogonal to the Chern characters of the matrix factorizations $k(m)^{\st}$, for $m\in\Z$.
Here $k(m)^{\st}$ is the stabilization of the trivial module $k$, with the grading shifted by $m$.
This is used in proving the analog of the Hodge-Riemann relations for matrix factorizations (see Proposition \ref{HR-prop}), 
which states that for a nonzero class $x\in H(W)^{\mu_d}\cap\HH_j(W)$ one has $\lan x,\ov{x}\ran_W\neq 0$
(and in fact, $\lan x,\ov{x}\ran_W$ is a positive multiple of a certain power of $i$). Roughly speaking, this is proved 
by reducing to the classical Hodge-Riemann relations for cohomology classes on the projective hypersurface $X$,
using Orlov's theorem relating the categories $\DMF_{\G_m}(W)$ and $D^b(X)$. For example, in the case $d=n$ Orlov's theorem
states an equivalence of these categories, and we check that $H(W)^{\mu_d}$ corresponds precisely to the primitive cohomology of $X$ under the corresponding isomorphism between the Hochschild homology (see Remark \ref{CY-rem}).

Since the image of $\psi$ is contained in $H(W)^{\mu_d}$, the above Hodge-Riemann relations show that the
vanishing \eqref{psi-vanishing-eq} is equivalent to the vanishing
$$\lan\ov{\psi(y')},\psi(y')\ran_W=0,$$
where $y'=\kappa(\a\cdot\Td_X^{-1})\cdot y$ and $y\in H^{p,q}$ with $p+q\neq 2\dim X-D$.
Using the adjointness of $\phi$ and $\psi$ we rewrite this as a certain purity property for the
composition $\phi\psi=(\Phi\circ\Psi)_*$. 

Finally, a computation with the Fourier-Mukai kernels (see Sections \ref{mf-ker-sec} and \ref{Koszul-mf-calc-sec})
shows that the kernel $K$ on $X\times X$ defining the functor $\Phi\circ\Psi: D^b(X)\to D^b(X)$ is given by the push-forward
of a Koszul matrix factorization of zero on some vector bundle over $X\times X$, supported along the zero section.
Applying the results of \cite{PV-Wc} (see Prop.\ \ref{homog-class-prop}) we obtain that the appropriate twist of the class
$\ch(K)$ is pure of certain dimension (see \eqref{main-homog-eq}). This gives the required purity property for $\phi\psi$
and so finishes the proof.

We conjecture that a statement similar to Theorem \ref{main-thm} holds for any quasihomogeneous $W$ 
with an isolated singularity. One could try to mimic our proof in the homogeneous case.
However, at present several technical ingredients are lacking. For example, in this case the analog of the
projective hypersurface is a DM-stack, so we need to identify the canonical bilinear form on the Hochschild homology
of a smooth proper DM-stack in terms of the Hochschild-Kostant-Rosenberg isomorphism. Another problem
is matching the effect of Orlov's equivalence on the Hochschild homology with the ad hoc isomorphism constructed in
\cite{CR}.

The paper is organized as follows. In Section \ref{prelim-sec} we build the background and prove some technical statements. 
The most important bits are
Section \ref{mf-ker-sec}, where we establish an adjunction result for Fourier-Mukai type functors involving matrix factorizations,
 and Section \ref{Orlov-sec} containing some calculations with Orlov's
equivalence. Then in Section \ref{rat-lattice-sec} we discuss rational lattices and prove the analog of the Hodge-Riemann
bilinear relations for the Hochschild homology of the category of matrix factorizations of a homogeneous polynomial with
an isolated singularity (or a quasihomogeneous polynomial as in Theorem \ref{main-thm}; see Proposition \ref{HR-prop}). 
Finally, in Section \ref{homog-sec} we recall the purity result from \cite{PV-Wc} and
show how Theorem \ref{main-thm} is deduced from it. We then deduce Theorem \ref{dim-property-thm} in Section
\ref{dim-property-sec}.
In the Appendix we prove a technical result involving
Grothendieck duality and matrix factorizations, which is needed in Section \ref{mf-ker-sec}.
\medskip

\noindent
{\it Conventions}.  We work with schemes and dg-categories over a field $k$. 
Starting from Section \ref{rat-lattice-sec}
we assume that $k=\C$. For a smooth projective variety $X$ we denote by $D^b(X)$ the bounded derived category of
coherent sheaves on $X$, which we equip with one of the standard dg-enhancements. We denote by $\Per(X)\sub D^b(X)$
the full subcategory of perfect complexes. For an algebraic group $G$ acting on $X$ we denote by $\Per_G(X)$ the
category of $G$-equivariant perfect complexes.
For a morphism of schemes $f$ we denote by $f^*$ and $f_*$ the corresponding
derived functors of pull-back and push-forward. 
For an additive category $\CC$ we denote by $\ov{\CC}$ its Caroubian completion.
By $\ch(\cdot)$ and $\Td(\cdot)$ we denote the characteristic classes of algebraic vector bundles constructed using
the Atiyah class, as in \cite{Atiyah} (their components differ from the corresponding topological characteristic classes by
factors of $(-2\pi i)$).

\medskip

\noindent
{\it Acknowledgments}. I would like to thank Alessandro Chiodo and Dimitri Zvonkine for organizing the workshop
``Mirror symmetry and spin curves" which prompted me to write this paper. I am grateful to Nick Rosenblyum for a useful discussion on adjoint functors. I am also grateful to the Mathematical Sciences Research Institute at Berkeley, 
where some of this work was done, for the hospitality and excellent working conditions. Finally, I thank the anonimous referee
for useful suggestions. This research is supported in part by the NSF grant DMS-1400390.


\section{Preliminaries}\label{prelim-sec}

\subsection{DG-categories and dg-functors}\label{dg-sec}

Let $k$ be a field. For a dg-category $\CC$ over $k$ we denote by $D(\CC)$ the derived category of right $\CC$-modules
and by $\Per(\CC)\sub D(\CC)$ the subcategory of perfect modules. We denote by $\Per_{dg}(\CC)$ the natural
dg-enhancement of $\Per(\CC)$ (defined using cofibrant right $\CC$-modules, see \cite[Sec.\ 7]{Toen}).

Throughout this paper we consider only $\CC$ such that $\Per_{dg}(\CC)$ is saturated (see \cite[Sec.\ 2.2]{TV}).

For an object $K\in \Per(\CC^{op}\ot\DD)$ we have a dg-functor of tensoring with $K$,
$$\Phi_K:\Per_{dg}(\CC)\to \Per_{dg}(\DD):M\mapsto M\ot^{\dL}_{\CC} K.$$
It is known that in this way we get a bijection between isomorphism classes in $\Per(\CC^{op}\ot\DD)$ and
morphisms between $\Per_{dg}(\CC)$ and $\Per_{dg}(\DD)$ in the localized
category $\Hqe$ of dg-categories,  obtained by inverting quasi-equivalences (see \cite{Toen}).
To a usual dg-functor $\Phi:\CC\to\DD$ we associate a bimodule, i.e., an object in $D(\CC^{op}\ot\DD)$,
$$(C,D)\mapsto \Hom_\DD(D,\Phi(C)),$$
which is perfect under our assumptions on $\CC$ and $\DD$.
The corresponding tensor functor $\Per_{dg}(\CC)\to\Per_{dg}(\DD)$ is an extension of $\Phi$.

In the remainder of the paper we often switch between $K$ and $\Phi_K$, and sometimes, denote them
by the same letter. We denote by $\De_\CC\in\Per(\CC^{op}\ot\CC)$ the {\it diagonal bimodule} corresponding to
the identity functor. We also have the corresponding dg-functor 
\begin{equation}\label{Tr-eq}
\Tr_{\CC}^{dg}:\Per_{dg}(\CC^{op}\ot \CC)\to\Per_{dg}(k)
\end{equation}
mapping $A^\vee\ot B$ to $\Hom_\CC(A,B)$.

For a pair of kernels 
$K_1\in\Per(\CC_1^{op}\ot\DD_1)$, $K_2\in\Per(\CC_2^{op}\ot\DD_2)$
we have the induced functor
$$K_1\Box K_2:\Per_{dg}(\CC_1\ot\CC_2)\to\Per_{dg}(\DD_1\ot\DD_2),$$
given by the external tensor product of the kernels $K_1$ and $K_2$.
Recall that there is a natural equivalence
$$\Per(\CC)^{op}\rTo{\sim} \Per(\CC^{op}):M\mapsto M^\vee$$
(see e.g., \cite[(3.6)]{Shk}). Using this equivalence, for a kernel $F\in\Per(\CC^{op}\ot\DD)$ we can define 
$F^{op}\in\Per(\CC\ot\DD^{op})$ as the kernel corresponding to the functor
$$\Per(\CC^{op})\rTo{\sim}\Per(\CC)^{op}\rTo{F}\Per(\DD)^{op}\rTo{\sim}\Per(\DD^{op}).$$


\begin{defi}\label{adjoint-defi}
Let $F\in\Per(\CC^{op}\ot\DD)$, $G\in\Per(\DD^{op}\ot\CC)$.
We say that $(F,G)$ is an {\it adjoint pair}, or that $F$ {\it is left adjoint to} $G$ if a morphism 
$$\varphi:\De_\CC\to G\circ F$$
is given in $D(\CC^{op}\ot\CC)$, such that for any $C\in\Per(\CC)$ and $D\in\Per(\DD)$ the induced morphism 
\begin{equation}\label{adj-eq1}
\Hom_\DD(F(C),D)\rTo{G} \Hom_\CC(GF(C),G(D))\rTo{?\circ\varphi} \Hom_\CC(C,G(D))
\end{equation}
is a quasi-isomorphism.
\end{defi}

Note that if $(F,G)$ is an adjoint pair then the map \eqref{adj-eq1} can be extended to a similar map of dg-functors
$\Per_{dg}(\CC^{op}\ot\DD)\to\Per_{dg}(k)$,
\begin{equation}\label{adj-eq}
\Tr_{\DD}\circ(F^{op}\Box\De_{\DD})\to \Tr_{\CC}\circ(\Id_{\CC}\Box G),
\end{equation}
induced by a quasi-isomorphism of kernels in $\Per(\CC\ot\DD^{op})$.
Indeed, the fact that \eqref{adj-eq} is a quasi-isomorphism of kernels is equivalent to the assertion that it becomes
a quasi-isomorphism of complexes when applied to $C^\vee\ot D\in \CC^{op}\ot\DD$, which corresponds to
\eqref{adj-eq1} being a quasi-isomorphism.

Recall that a {\it $\Z/2$-dg-category} is a dg-category in which all $\Hom$-complexes are $2$-periodic. 
Equivalently, we can replace the $2$-periodic complexes by the corresponding $\Z/2$-graded complexes.
The theory of $\Z/2$-dg-categories is parallel to the theory of dg-categories (see \cite[Sec.\ 4.1]{Dyck}).
With each dg-category $\CC$ we can associate its {\it $\Z/2$-folding} $\CC^{(2)}$, which is a $\Z/2$-dg-category such that
$$\Hom^{0}_{\CC^{(2)}}(A,B)=\bigoplus_{i\in\Z}\Hom^{2i}_{\CC}(A,B), \ \ \Hom^{1}_{\CC^{(2)}}(A,B)=\bigoplus_{i\in\Z}
\Hom^{2i+1}_{\CC}(A,B).$$
One can think of this operation as the tensor product with the algebra $k[u,u^{-1}]$, where $\deg(u)=2$.
It is easy to see that the Hochschild homology functor commutes with passing to the $\Z/2$-folding, i.e.,
the Hochschild homology of $\CC^{(2)}$ as a $\Z/2$-dg-category is the $\Z/2$-folding of the complex $HH_*(\CC)$,
and these identifications are compatible with the maps induced by functors.
In particular, taking the Chern character and the canonical pairing on the Hochschild homology discussed below are also
compatible with this operation.


\subsection{Chern character in Hochschild homology versus topological Chern character}\label{Chern-sec}

Recall that for any dg-category $\CC$ over $k$ and an object $E$ of $\CC$ one has 
the {\it Chern character} $\ch(E)\in HH_0(\CC)$ defined by the functoriality of the Hochschild homology (see e.g., \cite{Shk},
\cite{PV-HRR}).

In the case when $\CC=D^b(X)$, the dg-version of the bounded
derived category of coherent sheaves on a smooth projective
variety $X$ over $\C$, we have the Hochschild-Kostant-Rosenberg isomorphism
\begin{equation}\label{HKR-isom}
HH_*(D^b(X))=HH_*(X)\simeq \oplus_{p,q} H^q(X,\Om^p)
\end{equation}
so that $HH_i$ corresponds to the sum of terms with $p-q=i$.
Hence, by Hodge theory, we can identify the Hochschild homology $HH_*(X)$
with $H^*(X,\C)$.



As was shown in \cite{Cal}, 
under this identification, the abstract Chern character with values in $HH_*(X)$
is essentially the same
as the topological Chern character with values in $H^*(X,\C)$ (for $k=\C$).
Here is a more precise statement.

\begin{prop}\label{Chern-prop} Let $k=\C$.
For $E\in D^b(X)$ let $\ch^{top}(E)\in H^*(X,\C)$ denote the usual topological Chern character.
Then one has
$$\ch_k(E)=(-2\pi i)^k\ch_k^{top}(E).$$
\end{prop}

\Pf . Caldararu computes $\ch(E)$ in terms of traces of powers of the Atiyah class $\At(E)\in \Ext^1(E,E\ot\Om^1)$
(see the proof of \cite[Thm.\ 4.5]{Cal}). Thus, the assertion follows from the comparison 
with the topological Chern character in \cite[Sec.\ 5]{Atiyah}.\footnote{The factor $(-2\pi i)^k$ appears due to the standard
normalization of $\ch^{top}(E)$, so that it takes values in $H^*(X,\Z)$. Caldararu uses a different normalization, 
so he does not have this factor.}
\ed

In what follows we always use the notation $\ch(\cdot)$ for the abstract Chern character defined using Hochschild functoriality
(or in terms of the Atiyah class, in the case of $D^b(X)$). We denote by $\Td(\cdot)$ the Todd class defined by the
standard formulas in terms of the components of $\ch(\cdot)$ (we only use it for $D^b(X)$).

\subsection{Canonical pairing on the Hochschild homology}\label{can-pair-sec}


Applying the functoriality of the Hochschild homology to the functor \eqref{Tr-eq} and using 
the K\"unneth isomorphism for Hochschild homology we obtain a canonical pairing
\begin{equation}\label{pair-eq}
\lan\cdot,\cdot\ran_{\CC}:HH_*(\CC^{op})\ot HH_*(\CC)\to k
\end{equation}
(cf. \cite[Sec.\ 1.2]{Shk}).

Note that there is a canonical isomorphism $HH_*(\CC^{op})\simeq HH_*(\CC)$, so we can think of the canonical
pairing as a pairing on $HH_*(\CC)$.
However, one should be careful that in the case when the category $\CC$ is equipped with the duality $\CC^{op}\simeq\CC$
the induced identification of $HH_*(\CC^{op})\simeq HH_*(\CC)$ may be different from the canonical one.
For example this is the case for $\CC=D^b(X)$, where it is customary to use the duality on sheaves to identify $\CC$ with
$\CC^{op}$. This is related to the involution $a\mapsto \kappa(a)$ on the Hochschild homology of $X$ arising below
(cf. \cite[Rem.\ 2.15]{P-Lef}).

In the case of $\CC=D^b(X)$ (where $X$ is a smooth projective variety) 
the canonical pairing on $HH_*(D^b(X))$ 
is given in terms of the HKR-isomorphism \eqref{HKR-isom} by the formula
\begin{equation}\label{Ram-formula}
\lan a,b\ran_{D^b(X)}=(\kappa(a),b)_X, 
\end{equation}
where 
\begin{equation}\label{X-pairing}
(a,b)_X:=\int_X a\cdot b\cdot\Td_X,
\end{equation}
\begin{equation}\label{iota-eq}
\kappa(c)=(-1)^q c \text{ for } c\in H^{p,q}=H^q(X,\Om^p)
\end{equation}
(see \cite[Eq.\ (8)]{Ram}; note that Ramadoss in \cite{Ram} works with the Mukai pairing which coincides
with the canonical pairing \eqref{pair-eq} for $\CC=D^b(X)$, e.g., by \cite[Prop.\ 2.14]{P-Lef}).
Here $\int_X$ is defined as the projection
$$\int_X:\oplus_{p,q} H^q(X,\Om^p)\to H^n(X,\Om^n_X)\simeq k,$$
where $n=\dim X$. 
In the case $k=\C$ it is related to the topological operation of integration over the fundamental cycle of $X$
by the formula
$$\int_X\om=\frac{1}{(2\pi i)^n}\int^{top}_{[X]}\om$$
(see e.g., \cite{Del-HC}).

The classical Hodge-Riemann bilinear relations (see e.g., \cite[Ch.\ V]{Wells})
imply the following property of the canonical pairing $\lan\cdot,\cdot\ran_{D^b(X)}$.

\begin{lem}\label{HR-lem} 
Let $H\in\Pic(X)$ be an ample class, and let
$P^n_H(X)\sub H^n(X)$ be the corresponding primitive part of the middle cohomology consisting
of the classes $a$ such that $a\cdot c_1(H)=0$.
Then for a nonzero class $a\in P^n_H(X)$
one has
$$(-1)^{\frac{n(n+1)}{2}}\lan a,\ov{a}\ran_{D^b(X)}>0.$$
\end{lem}

\Pf . By \eqref{Ram-formula}, we have
$$\lan a,\ov{a}\ran_{D^b(X)}=(-1)^q\int_X a\cdot\ov{a}\cdot\Td_X=\frac{(-1)^q}{(2\pi i)^n}\int^{top}_{[X]}a\cdot\ov{a}\cdot\Td_X,$$
where $a\in H^{p,q}$, $p+q=n$. Since in our case $a\cdot\ov{a}\in H^{n,n}(X)$, we can delete $\Td_X$,
and the result follows from the Hodge-Riemann bilinear relations stating that
$$(-1)^{\frac{n(n-1)}{2}}i^{p-q}\int^{top}_{[X]}a\cdot\ov{a}>0.$$
\ed


\subsection{Maps on Hochschild homology induced by Fourier-Mukai functors}\label{Hoch-maps-sec}




We need the fact that adjoint functors induce adjoint operators on Hochschild homology
with respect to the canonical pairings \eqref{pair-eq} (cf.\  \cite[Thm.\ 8]{CW} in the case of derived categories of sheaves).

\begin{lem}\label{adjoint-HH-lem} 
If $(F\in\Per(\CC^{op}\ot\DD), G\in\Per(\DD^{op}\ot\CC))$ is an adjoint pair of kernels then
$$\lan F_*(x),y\ran_{\DD}=\lan x,G_*(y)\ran_{\CC},$$
where $x\in HH_*(\CC)$, $y\in HH_*(\DD)$, and
$F_*:HH_*(\CC)\to HH_*(\DD)$, $G_*:HH_*(\DD)\to HH_*(\CC)$ are the induced maps on the Hochschild
homology.
\end{lem}

\Pf . The required equality is immediately obtained from \eqref{adj-eq} by passing to the induced maps on
Hochschild homology. We have to use the fact that the map
$$F^{op}_*:HH_*(\CC^{op})\to HH_*(\DD^{op})$$
coincides with $F_*$ under the natural identifications $HH_*(\CC)\simeq HH_*(\CC^{op})$, 
$HH_*(\DD)\simeq HH_*(\DD^{op})$. The simplest way to check this is to use the definition of the maps $F_*$ via
the Hochschild chain complexes.
\ed

We will use the following formula for the
maps induced on the Hochschild homology by Fourier-Mukai type functors,
in terms of the Chern character of the kernel and the pairing \eqref{X-pairing}.

\begin{lem}\label{FM-HH-lem} 
For the Fourier-Mukai functor $F:D^b(X)\to D^b(Y)$ associated with a kernel
$K\in D^b(X\times Y)$ the induced map on the Hochschild homology gets identified via the HKR-isomorphisms
with 
$$F_*:H^*(X,\C)\to H^*(Y,\C): a\mapsto \tr_{12}(a\ot \ch(K)),$$
where 
$$\tr_{12}:H^*(X)\ot H^*(X)\ot H^*(Y)\to H^*(Y): a\ot b\ot c\mapsto (a,b)_X c.$$
\end{lem}

\Pf . This is equivalent to \cite[Thm.\ 1.2]{MS}.
\ed

\subsection{Matrix factorizations}\label{mf-sec}

In this section we recall some basic results about matrix factorizations and also prove several technical statements that
will be needed later in working with Fourier-Mukai transforms involving matrix factorizations. The key result that allows
to deduce many results for matrix factorizations from the classical results about coherent sheaves is the equivalence with
the singularity category of the hypersurface $W=0$ (see \eqref{singularity-category-equiv} below).

Let $X$ be a $\G_m$-scheme, i.e., a scheme with a $\G_m$-action.
Throughout this paper we make an assumption that our $\G_m$-schemes admit a $\G_m$-equivariant ample line bundle. 
This implies that they admit a $\G_m$-invariant open affine cover and have a resolution property, i.e., 
every $\G_m$-equivariant coherent sheaf on $X$ admits a surjection from
a $\G_m$-equivariant vector bundle (see \cite{Thom}).

For a $\G_m$-equivariant quasicoherent sheaf $\FF$ on $X$ we denote
by $\FF\{i\}$ the same sheaf with the $\G_m$-action twisted by the character $\la\mapsto \la^i$ of $\G_m$.

Let $W$ be a function on $X$ of weight $d>0$ with respect to the $\G_m$-action, i.e., $W\in H^0(X,\OO_X\{d\})^{\G_m}$.
A {\it $\G_m$-equivariant matrix factorization of $W$} is a $\Z/2$-graded vector bundle $E=E_0\oplus E_1$ on $X$
together with maps 
$$\de_1:E_1\to E_0, \ \ \de_0:E_0\to E_1\{d\},$$
such that $\de_0\de_1=W\cdot\id$ and $\de_1\de_0=W\cdot\id$.
For a pair of $\G_m$-equivariant matrix factorizations of $W$, $E$ and $F$ we consider the complex of
$\G_m$-equivariant sheaves $\und{\Hom}(E,F)$ given by
$$\und{\Hom}^{2n}(E,F)=\und{\Hom}(E_0,F_0\{dn\})\oplus\und{\Hom}(E_1,F_1\{dn\}),$$
$$\und{\Hom}^{2n+1}(E,F)=\und{\Hom}(E_0,F_1\{d(n+1)\})\oplus\und{\Hom}(E_1,F_0\{dn\}),$$
with the differential $f\mapsto \de_F\circ f-(-1)^{|f|}f\circ\de_E$.
Let 
$$R\Ga(X,\cdot):\Com(\Qcoh_{\G_m}(X))\to \Com(k[\G_m])$$
be a multiplicative dg-model of the push-forward to the point, given e.g., by the Cech complex
with respect to a $\G_m$-invariant open affine cover (as in \cite[Sec.\ 2]{Shipman}).
We define the dg-category $\MF_{\G_m}(X,W)$ of $\G_m$-equivariant matrix factorizations of $W$ by setting
$$\Hom(E,F)=R\Ga(X,\und{\Hom}(E,F))^{\G_m}.$$
Passing to the $0$th cohomology of the dg-category $\MF_{\G_m}(X,W)$
we get the derived category of matrix factorizations $\DMF_{\G_m}(X,W)$.

Assume that $X$ is smooth. Recall that by our assumption $X$ admits a $\G_m$-equivariant ample line bundle,
hence it has a $\G_m$-resolution property. Assume also that $W$ is not a zero divisor.
Then the functor associating with a matrix factorization $E$ the cokernel of $\de_1:E_1\to E_0$ extends to an equivalence
\begin{equation}\label{singularity-category-equiv}
\DMF_{\G_m}(X,W)\rTo{\sim} D_{\Sg,\G_m}(X_0),
\end{equation}
where $X_0\sub X$ is the hypersurface $W=0$, $D_{\Sg,\G_m}(X_0)$ is the {\it singularity category}, defined as the quotient
of the bounded $\G_m$-equivariant derived category by the subcategory $\Per_{\G_m}(X_0)$ of perfect complexes.
In this form the equivalence follows from \cite[Thm.\ 3.14]{PV-mf-stack} but the construction and the main ideas go back to Orlov \cite{Orlov-dbr} (see also \cite[Sec.\ 3]{Orlov-sing}, \cite{Orlov-nonaff} and \cite{EPos}).

Note that there is a different way to define the derived category of matrix factorizations as the {\it absolute derived
category} $D^{\abs}(\MF_{\G_m}(X,W))$, which is the quotient of the naive homotopy category by the convolutions
of exact sequences of matrix factorizations (see \cite{Orlov-nonaff}, \cite{EPos}). The equivalence of this definition
with the one above (in the case of smooth $X$) follows from the equivalence
$$D^{\abs}_{\G_m}(X,W)\simeq D_{\Sg,\G_m}(X_0)$$
(which can be proved as in \cite{Orlov-nonaff}, \cite{EPos}).
In the case when $X$ is affine (and smooth) the derived category of matrix factorization on $X$ coincides with the naive
homotopy category (see e.g., \cite[Sec.\ 3]{Orlov-sing}, \cite[Prop.\ 3.19]{PV-mf-stack}).

For technical reasons we often work with the Caroubian completion $\ov{\DMF}_{\G_m}(X,W)$ of the derived category
of matrix factorizations. It can be realized as a full subcategory in the derived category of matrix factorizations
of quasicoherent sheaves (see \cite[Sec.\ 4]{PV-mf-stack}, \cite[Sec.\ 2.3]{EPos}). In the case when $W$ is a homogeneous
polynomial on $\A^n$ with an isolated singularity the relation of $\DMF_{\G_m}(\A^n,W)$ with the
derived category $D^b(Y)$ on the corresponding projective hypersurface $Y$ (see Sec.\ \ref{Orlov-sec}) implies
that $\DMF_{\G_m}(\A^n,W)$ is in fact Caroubian closed.

We have a natural duality equivalence
$$\Du: \MF_{\G_m}(X,W)^{op}\to \MF_{\G_m}(X,-W),$$
where $\Du(E)_0=E_0^\vee$, $\Du(E)_1=E_1^\vee\{-d\}$ with the induced differential
(see \cite[Eq.\ (2.13)]{PV-HRR}). Often, we will simply write $E^\vee$ instead of $\Du(E)$.
Similarly to the $\Z/2$-graded case (see \cite[Lem.\ 3.9]{LP}) this duality can be interpreted in terms
of the singularity category.

\begin{lem}\label{duality-lem} Assume that $X$ is smooth and $W$ is not a zero divisor. 
Under the equivalences
of $\DMF_{\G_m}(X,W)$ and $\DMF_{\G_m}(X,-W)$ with $D^b_{\Sg,\G_m}(X_0)$, the duality $\Du$ corresponds
to the duality $F\mapsto R\und{\Hom}(F,\OO_{X_0}[-1])$.
\end{lem}

\Pf . Recall that for a matrix factorization $E$ the corresponding object of the singularity category is represented
by the coherent sheaf $F$ on $X_0$ fitting into the exact sequence
$$0\to E_1\to E_0\to i_*F\to 0,$$
where $i:X_0\to X$ is the embedding. By duality, we have an exact triangle
$$E_0^\vee\to E_1^\vee\to R\und{\Hom}(i_*F,\OO_X)[1]\to\ldots$$
By Grothendieck duality, 
$$R\und{\Hom}(i_*F,\OO_X)[1]\simeq i_*R\und{\Hom}(F,i^!\OO_X)[1]\simeq i_*R\und{\Hom}(F,\OO_{X_0}).$$
Thus, $G=R\und{\Hom}(F,\OO_{X_0})$ is a sheaf on $X_0$, and we have an exact sequence
$$0\to E_0^\vee\to E_1^\vee\to i_*G\to 0.$$
On the other hand, the object of the singularity category associated with $\Du(E)$ is the coherent sheaf $F'$ on $X_0$
from the exact sequence
$$0\to E_1^\vee\{-d\}\to E_0^\vee \to i_*F'\to 0.$$
Hence, we have an exact sequence on $X_0$
$$0\to F'\to i^*E_1^\vee \to G\to 0$$
which shows that $F'\simeq G[-1]$ in the singularity category.
\ed

For a pair of potentials $W, W'$ on $X$, both of weight $d>0$, we define
the tensor product functor
$$\otimes:\MF_{\G_m}(X,W)\times \MF_{\G_m}(X,W')\to \MF_{\G_m}(X,W+W')$$
by setting
\begin{equation}\label{usual-tensor-product}
(E\ot F)_0=E_0\ot F_0\oplus E_1\ot F_1\{d\} \ \text{ and } (E\ot F)_0=E_0\ot F_1\oplus E_1\ot F_0
\end{equation}
with the differential $\de_E\ot\id_F+\id_E\ot \de_F$.

For $F\in\MF_{\G_m}(X,W)$ we can consider the infinite complex of $\G_m$-equivariant sheaves on 
$X_0=W^{-1}(0)$,
$$\com(F): \ldots\to E_0\{-d\}|_{X_0}\to E_1|_{X_0}\to E_0|_{X_0}\to E_1\{d\}|_{X_0}\to \ldots$$
with $E_0|_{X_0}$ placed in degree $0$. Note that if $W=0$ then this is a complex on $X=X_0$.
The following relation between $\Hom$'s, duality and tensor product is straightforward to check
(cf.\ \cite[Lem.\ 1.1.6]{PV-CohFT}).

\begin{lem}\label{Hom-lem} 
For $E,F\in\MF_{\G_m}(X,W)$ one has an isomorphism of $\G_m$-equivariant complexes
on $X$,
$$\und{\Hom}(E,F)\simeq \com(E^\vee\ot F).$$
Hence,
$$\Hom(E,F)\simeq R\Ga(X,\com(E^\vee\ot F))^{\G_m}.$$
\ed
\end{lem}

For a $\G_m$-scheme $X$ and for some $d>0$ we can consider the category
$\MF_{\G_m,d}(X,0)$
of matrix factorizations of $0$ on $X$, where $0$ is viewed as a function of weight $d$.
In this situation we define the functor 
\begin{equation}\label{maf-eq}
\maf: \Per_{\G_m}(X)\to\DMF_{\G_m,d}(X,0): \ \ \maf(C^\bullet)_0=\bigoplus_n C^{2n}\{-nd\}, \ \ \maf(C^\bullet)_1=\bigoplus_n C^{2n-1}\{-nd\}
\end{equation}
(note that since we have a resolution property for $\G_m$-equivariant sheaves on $X$, an object of $\Per_{\G_m}(X)$ can be
represented globally by a bounded complex of $\G_m$-equivariant vector bundles, see \cite[Lem.\ 3.5]{PV-mf-stack}).
If $W$ is a function on $X$ of weight $d$ then we also get the tensor product operation 
$$\otimes:\Per_{\G_m}(X)\times \DMF_{\G_m}(X,W)\to \DMF_{\G_m}(X,W), \ \ F\ot E:=\maf(F)\ot E.$$
Then one has a natural isomorphism of complexes on $X_0$,
\begin{equation}\label{com-tensor-eq}
\com(F\ot E)\simeq F|_{X_0}\ot\com(E)
\end{equation}
(see \cite[Lem.\ 1.1.5]{PV-CohFT}). This easily implies (assuming $W$ is not a zero divisor) that
under the equivalence \eqref{singularity-category-equiv} the operation $F\ot ?$ corresponds to
the operation $F|_{X_0}\ot ?$ on the singularity category $D_{\Sg,\G_m}(X_0)$.


Assume now in addition that the action of $\G_m$ on $X$ is trivial and $X$ is quasiprojective.
Then we can define several complexes associated with a $\G_m$-equivariant matrix factorization
of $0$ on $X$. 
Note that for any $\G_m$-equivariant matrix factorization $E$ of $0$
we can write $E_0=\bigoplus E_{0,i}$, $E_1=\bigoplus E_{1,i}$, where $\G_m$ acts on $E_{\bullet,i}$ through the character 
$\la\mapsto \la^{-i}$, i.e., $E_{\bullet,i}=(E_\bullet\{i\})^{\G_m}$. Then we have the functor
$$\com_0:\DMF_{\G_m,d}(X,0)\to \Per(X),$$
where $\com_0(E)$ is the complex
$$\ldots\to E_{1,0}\to E_{0,0}\to E_{1,d}\to\ldots$$
with $E_{0,0}$ placed in degree $0$.
Note that $E_{0,i}$ and $E_{1,i}$ are nonzero only for finitely many $i$, so the complex $\com_0(E)$ is bounded.
Let us also consider the functor 
\begin{equation}\label{ov-com-eq}
\ov{\com}:\DMF_{\G_m,d}(X,0)\to \Per(X), \ \ \ov{\com}(E)=\bigoplus_{i=0}^{d-1} \com_0(E\{i\}).
\end{equation}

\begin{prop}\label{com-maf-prop} 
Let $X$ be a scheme equipped with the trivial $\G_m$-action. 

\noindent
(i) For $E\in\DMF_{\G_m,d}(X,0)$ one has a natural isomorphism of complexes 
$$\com_0(E)\simeq \com(E)^{\G_m}.$$

\noindent
(ii) The functor $\maf$ is left adjoint to $\com_0$.

\noindent
(iii) The left adjoint to the functor $\ov{\com}$ is $F\mapsto \bigoplus_{i=0}^{d-1}\maf(F)\{-i\}$.

\noindent
(iv) For $F\in \Per(X)$ and $E\in \DMF_{\G_m,d}(X,0)$ one has a natural isomorphism of complexes
$$\ov{\com}(F\ot E)\simeq F\ot\ov{\com}(E),$$
where on the left we equip $F$ with the trivial $\G_m$-action.
\end{prop}

\Pf . (i) This follows immediately from the definitions.

\noindent
(ii) For a complex $C^\bullet\in\Per(X)$ and a matrix factorization $F$ let us compute the complex
$\und{\Hom}(\maf(C^\bullet),F)^{\G_m}$. We have
\begin{align*}
&\und{\Hom}^{2n}(\maf(C^\bullet),F)^{\G_m}=
\und{\Hom}(\bigoplus_i C^{2i}\{-id\},F_0\{dn\})^{\G_m}\oplus
\und{\Hom}(\bigoplus_n C^{2i-1}\{-id\},F_1\{dn\})^{\G_m}=\\
&\prod_i\und{\Hom}(C^{2i},F_{0,d(i+n)})\oplus\prod_i\und{\Hom}(C^{2i-1}F_{1,d(i+n)})=
\prod_j\und{\Hom}(C^j,\com_0(F)^{j+2n}).
\end{align*}
Similarly
$$\und{\Hom}^{2n+1}(\maf(C^\bullet),F)^{\G_m}=\prod_j\und{\Hom}(C_j,\com_0(F)^{j+2n+1}),$$
and the differentials match. Passing to $R\Ga(X,?)$ we get the required adjointness.

\noindent
(iii) This follows easily from (ii).

\noindent
(iv) This follows from \eqref{com-tensor-eq}: first, we check a similar property for $\com_0$, and then for $\ov{\com}$.
\ed

\begin{rem} When the action of $\G_m$ on $X$ is trivial, we have an equivalence
$$\DMF_{\G_m,d}(X,0)\simeq \Per_{\mu_d}(X)$$
associating to $E$ the complex $\bigoplus_{i=0}^{d-1}\com_0(E\{i\})\{-i\}$ (cf.\ \cite[Prop.\ 1.2.2]{PV-CohFT}).
The functor $\ov{\com}$ is the composition of this equivalence with the forgetful functor $\Per_{\mu_d}(X)\to \Per(X)$.
\end{rem}

Let us now return to the general situation of a $\G_m$-scheme $X$ with a function $W$ of weight $d>0$.
For a closed $\G_m$-invariant subset $T\sub X_0=W^{-1}(0)$ we denote by $\DMF_{\G_m,T}(X,W)\sub\DMF_{\G_m}(X,W)$
the full subcategory of matrix factorizations $E$ such that for every closed point $x\in X_0\setminus T$
the complex $\com(E)|_x$ is exact. Equivalently, these are matrix factorizations that become trivial in the category
$\DMF_{\G_m}(X\setminus T,W)$ (see \cite[Lem.\ 5.4(iii)]{PV-mf-stack}).

Let $f:X\to Y$ be a $\G_m$-morphism of smooth $\G_m$-varieties, and let $W$ be a function of weight $d>0$ on $Y$.
The push-forward functor $f_*$ for matrix factorizations is most naturally defined in terms of matrix factorizations of quasicoherent sheaves (see \cite[Sec.\ 3]{EPos}, \cite{BDFIK}). If $T\sub X$ is a closed subset, proper over $Y$, then one has also a natural push-forward functor
$$f_*:\ov{\DMF}_{\G_m,T}(X,f^*W)\to \ov{\DMF}_{\G_m}(Y,W)$$
(see \cite[Sec.\ 3.4]{EPos}, \cite[Sec.\ 6]{PV-mf-stack}). 

Assume in addition that $W$ and $f^*W$ are not zero divisors.
Let $f_0:X_0=f^{-1}(W^{-1}(0))\to Y_0=W^{-1}(0)$ be the morphism between the hypersurfaces induced by $f$.
Then the functor $f_{0*}:D^b(X_0)\to D^b(Y_0)$ induces a functor between the singularity categories, 
which corresponds to $f_*$ under the equivalences \eqref{singularity-category-equiv} for $W$ and $f^*(W)$.


The following property is straightforward.

\begin{lem}\label{com-f*-lem} 
For $f:X\to Y$ a $\G_m$-equivariant morphism of smooth $\G_m$-varieties, and for $F$, a $\G_m$-equivariant matrix factorization of $0$ on $X$, one has 
$$\com(f_*(F))\simeq f_*\com(F).$$
\end{lem}

The analog of the Grothendieck duality for matrix factorizations was established under
quite general assumptions in the work of Efimov and Positselski \cite[Sec.\ 3]{EPos}. Here we will use the following version for smooth morphisms between smooth varieties.

\begin{prop}\label{GD-prop} 
Let $f:X\to Y$ be a smooth $\G_m$-equivariant morphism of relative dimension $m$ between smooth $\G_m$-varieties,
and let $T\sub X$ be a closed subset, proper over $Y$. Let $W$ be a function on $Y$ of weight $d>0$ with
respect to the $\G_m$-action, which is not a zero divisor. 
For $E\in\ov{\DMF}_{\G_m,T}(X,f^*W)$ we have a natural functorial isomorphism
$$f_*(E^\vee\ot\DD_f)\simeq (f_*E)^\vee$$
in $\ov{\DMF}_{\G_m}(Y,-W)$, where
$$\DD_f:=\om_f[m].$$
\end{prop}

\Pf . Let $f_0:X_0\to Y_0=W^{-1}(0)$ be the morphism induced by $f$, where $X_0=f^{-1}(Y_0)$.
Note that under the equivalence \eqref{singularity-category-equiv}, extended to Caroubian completions,
$E$ corresponds to an object $F\in\ov{D}_{\Sg}(X_0)$ supported at $T\cap X_0$
(see \cite[Prop.\ 5.6]{PV-mf-stack}). Hence, by Lemma \ref{duality-lem}, it suffices to construct an isomorphism
$$f_{0*}(R\und{\Hom}(F,\OO_{X_0})\ot \DD_f|_{X_0})\simeq R\und{\Hom}(f_{0*}F,\OO_{Y_0})$$
in $D^b(Y_0)$, provided $F\in D^b(X_0)$ has a support proper over $Y_0$. Since $\DD_f|_{X_0}\simeq \DD_{f_0}$,
such an isomorphism is given by the usual Grothendieck duality.
\ed

\begin{cor}\label{duality-cor} 
In the situation of Proposition \ref{GD-prop}, for $E\in\ov{\DMF}_{\G_m,T}(X,f^*W)$ and $F\in\ov{\DMF}_{\G_m}(Y,W)$ one has
$$\und{\Hom}(f_*E,F)\simeq f_*\und{\Hom}(E,f^+F), \text{ where}$$
$$f^+F=\DD_f\ot f^*F.$$
Hence, we also have a functorial isomorphism
$$\Hom(f_*E,F)\simeq \Hom(E,f^+F).$$
\end{cor}

\Pf . Using Lemma \ref{Hom-lem}, Proposition \ref{GD-prop}, and the projection formula we get
$$\und{\Hom}(f_*E,F)\simeq \com((f_*E)^\vee\ot F)\simeq \com(f_*(E^\vee\ot \DD_f\ot f^*F)).$$
Using Lemma \ref{com-f*-lem} we can switch $\com$ with $f_*$, so we obtain
$$\und{\Hom}(f_*E,F)\simeq f_*\com(E^\vee\ot \DD_f\ot f^*F)\simeq f_*\und{\Hom}(E,\DD_f\ot f^*F),$$
as claimed. The second isomorphism is obtained from the first by applying $R\Ga$.
\ed 

Let us recall an important construction of {\it Koszul matrix factorizations}, which can be viewed as
a generalization of the Koszul complex. Assume we have a $\G_m$-equivariant vector bundle $V$ on $X$ and
invariant global sections 
$$\a\in H^0(X,V\{d\})^{\G_m},\ \b\in H^0(X,V^{\vee})^{\G_m} \ \text{ such that }\lan\a,\b\ran=W.$$ 
Then we define the {\it Koszul matrix factorization} 
$\{\a,\b\}$ of $W$ by
$$\{\a,\b\}_0=\OO_X\oplus \we^2 V\{d\}\oplus \we^4V\{2d\}\oplus\ldots,$$
$$\{\a,\b\}_1=V\oplus \we^3 V\{d\}\oplus \we^5V\{2d\}\oplus\ldots,$$
with the differential given by
\begin{equation}\label{Koszul-diff-eq}
\de_{\a,\b}=\a\we?+\iota_{\b},
\end{equation}
where $\iota_{\b}$ is the contraction by $\b$.
An important fact is that $\{\a,\b\}$ is supported on the locus of common zeros of $\a$ and $\b$ (see \cite[Lem.\ 1.5.1]{PV-CohFT}).



\subsection{Some functors given by kernels and an adjunction between them}\label{mf-ker-sec}

Suppose we have a diagram of smooth $\G_m$-varieties
\begin{equation}\label{ker-diag-eq}
\begin{diagram}
&&Y\\
&\ldTo{f}&&\rdTo{p}\\
Z&&&& X
\end{diagram}
\end{equation}
where $f$ is a smooth morphism of constant relative dimension, and $\G_m$ acts trivially on $X$.
Assume further that $W$ is a function on $Z$
of weight $d>0$ with respect to the $\G_m$-action. 
Given a $\G_m$-equivariant matrix factorization $P$ of $-f^*W$ on $Y$,
with proper support, we can define functors
\begin{equation}\label{Phi-P-eq}
\begin{array}{l}
\wt{\Phi}_{P}:\ov{\DMF}_{\G_m}(Z,W)\to \ov{\DMF}_{\G_m,d}(X,0), \ \ \Phi_{P}: \ov{\DMF}_{\G_m}(Z,W)\to D^b(X),\\
\wt{\Psi}_{P}:\ov{\DMF}_{\G_m,d}(X,0)\to \ov{\DMF}_{\G_m}(Z,-W), \ \Psi_P:D^b(X)\to  \ov{\DMF}_{\G_m}(Z,-W),
\end{array}
\end{equation}
as follows: 
$$\wt{\Phi}_{P}(E)=p_*(f^*E\ot P), \ \Phi_P=\ov{\com}\circ\wt{\Phi}_P,$$
$$\wt{\Psi}_{P}(F)=f_*(P\ot p^*F), \ \ \Psi_P(F)=\bigoplus_{i=0}^{d-1}\wt{\Psi}_{P}(\maf(F)\{-i\}).$$
Here $f^*E\ot P$ (resp., $P\ot p^*F$) is a $\G_m$-equivariant matrix factorization of $0$ (resp., $-f^*W$) on $Y$, that has a proper support, and so we can apply the push-forward functor $p_*$ (resp., $f_*$) to it.
The functors $\maf$ and $\ov{\com}$ are given by \eqref{maf-eq} and \eqref{ov-com-eq}.

\begin{lem}\label{gen-comp-lem} Let 
\begin{diagram}
&&Y'\\
&\ldTo{g}&&\rdTo{q}\\
Z&&&& X'
\end{diagram}
be another diagram with the same properties as \eqref{ker-diag-eq}, and let $P'$ be a $\G_m$-equivariant matrix
factorization of $g^*W$ on $Y'$ with proper support.
Then the composition 
$$\ov{\DMF}_{\G_m,d}(X',0)\rTo{\wt{\Psi}_{P'}} \ov{\DMF}_{\G_m}(Z,W)\rTo{\wt{\Phi}_P} \ov{\DMF}_{\G_m,d}(X,0)$$
is isomorphic to the Fourier-Mukai type functor associated with the kernel
$$\wt{K}=p_{X'X,*}(p_1^*P'\otimes p_2^*P)\in \ov{\DMF}_{\G_m,d}(X'\times X,0),$$ 
where we consider the diagram
\begin{diagram}
&&&&Y'\times_{Z} Y\\
&&&\ldTo{p_1}&&\rdTo{p_2}\\
&&Y'&&&&Y\\
&\ldTo{q}&&\rdTo{g}&&\ldTo{f}&&\rdTo{p}\\
X'&&&&Z&&&&X
\end{diagram}
and denote by $p_{X'X}:Y'\times_{Z}Y\to X'\times X$  
the map induced by $q\circ p_1$ and $p\circ p_2$.
The composition 
$$\Phi_P\circ\Psi_{P'}:D^b(X')\to D^b(X)$$
is isomorphic to the Fourier-Mukai functor associated with a kernel $K\in D^b(X'\times X)$ such that
$[K]=d[\ov{\com}(\wt{K})]$ in $K_0(X'\times X)$.
\end{lem}

\Pf . The first assertion follows easily from the projection formula and the base change formula (cf. \cite[Sec.\ 5.2]{BDFIK}). It remains to compute
the composition 
$$\Phi_P\circ\Psi_{P'}: F\mapsto \ov{\com}(\wt{\Phi}_P\circ\wt{\Psi}_{P'})(\bigoplus_{i=0}^{d-1} \maf(F)\{-i\}).$$
Let $p_X:X'\times X\to X$, $p_{X'}:X'\times X\to X'$ be the projections. We have
$$\ov{\com}\bigl(p_{X,*}(p_{X'}^*\maf(F)\{-i\}\ot \wt{K})\bigr)\simeq p_{X,*}\ov{\com}(p_{X'}^*F\ot \wt{K}\{-i\})\simeq
p_{X,*}\bigl(p_{X'}^*F\ot\ov{\com}(\wt{K}\{-i\})\bigr),$$
where in the last isomorphism we used Proposition \ref{com-maf-prop}(iv).
Thus, $\Phi_P\circ\Psi_{P'}$ is associated with the kernel
$$K=\bigoplus_{i=0}^{d-1}\ov{\com}(\wt{K}\{-i\}).$$
It is easy to see that $\ov{\com}(\wt{K}\{-i\})$ and $\ov{\com}(\wt{K})$ 
are direct sums of almost the same complexes---some complexes
get shifted by even integers. Hence, the assertion about the classes in the Grothendieck group.
\ed

As before, we denote $\DD_f=\om_f[\dim Y-\dim Z]$.

\begin{prop}\label{mf-adjoint-prop} 
In the above situation the left adjoint functor to $\wt{\Phi}_P$ (resp., $\Phi_P$) is $\wt{\Psi}_{Q}$ (resp., $\Psi_{Q}$), where 
$Q=P^\vee\ot \DD_f\in \MF_{\G_m}(Y,f^*W)$. Furthermore, the corresponding dg-functors are adjoint in the sense
of Definition \ref{adjoint-defi}.
\end{prop}

\Pf . For $E\in \MF_{\G_m,d}(X,0)$ and $F\in\MF_{\G_m}(Z,W)$,
we have a chain of quasi-isomorphisms 
\begin{align*}
&\Hom(E,\wt{\Phi}_P(F))\simeq 
\Hom(E,p_*(P\ot f^*F))\simeq \\
&\Hom(p^*E,P\ot f^*F)\simeq\Hom(p^*E\ot P^\vee\ot\DD_f, f^*F\ot \DD_f).
\end{align*}
Note that $p^*E\ot P^\vee\ot\DD_f$ has proper support.
Hence, using Corollary \ref{duality-cor}  we can rewrite this as
$$\Hom(f_*(p^*E\ot P^\vee\ot\DD_f),F)\simeq \Hom(\wt{\Psi}_{Q}(E),F).$$
To show the required adjunction at the dg-level we have to prove that 
the above isomorphism is induced by a map $\Id\to \wt{\Phi}_P\circ\wt{\Psi}_Q$ given by a map of kernels.
By Lemma \ref{gen-comp-lem}, the composition $\wt{\Phi}_P\circ\wt{\Psi}_Q$ is the functor given by the kernel
$$p_{XX,*}(p_1^*Q\ot p_2^*P)\in  \MF_{\G_m,d}(X'\times X,0),$$ 
where $p_1$ and $p_2$ are the two projections $Y\times_{Z} Y\to Y$, and
$$p_{XX}:Y\times_{Z} Y\to X\times X$$
is the map with the components $(pp_1,pp_2)$.
Thus, we need a map of kernels
$$\De_*\OO_X\to p_{XX,*}(p_1^*Q\ot p_2^*P).$$
The construction of this map and the verification 
that the corresponding natural transformation $\Id_{\MF_{\G_m,d}(X,0)}\to \wt{\Phi}_P\circ\wt{\Psi}_Q$ induces
the same quasi-isomorphism
$$\Hom(\wt{\Psi}_{Q}(E),F)\rTo{\sim}\Hom(E,\wt{\Phi}_P(F))$$
as the one obtained above are done in the Appendix (see \eqref{map-of-kernels} and Proposition \ref{compatibility-prop}).

The adjunction of the pair $(\Psi_{Q},\Phi_P)$ follows from the adjunction of the pair $(\wt{\Psi}_Q,\wt{\Phi}_P)$
using Proposition \ref{com-maf-prop}(iii).
\ed

Now assume that in the above situation we have $Z=\A^n$ and $W$ is a quasihomogeneous polynomial on $Z$
of degree $d>0$ with an isolated singularity.
Recall that we have a decomposition 
\begin{equation}\label{HH-W-decomp}
HH_*(\MF_{\G_m}(W))=\HH(W)=\bigoplus_{\ga\in \mu_d} H(W_\ga)^{\mu_d},
\end{equation}
where $W_\ga=W|_{(\A^n)^\ga}$ (see \cite[Thm.\ 2.6.1]{PV-CohFT}).
In fact, this decomposition is exactly the $\mu_d$-grading associated with the natural
$\Z/d$-action on $HH_*(\MF_{\G_m}(W))$ induced by the functors $E\mapsto E\{i\}$ (see \cite[Thm.\ 2.6.1(ii)]{PV-CohFT}). 
Let $\Pi$ denote the projector of $\HH(W)$  onto the summand $H(W)^{\mu_d}$, corresponding to $\ga=1$.

\begin{lem}\label{phi-inv-lem}
Let $\phi_P:\HH(W)\to HH_*(X)=H^*(X,\C)$ be the map induced on Hochschild homology by the functor $\Phi_P$
given by \eqref{Phi-P-eq}.
We have $\phi_P=\phi_P\circ \Pi$.
\end{lem}

\Pf . It is clear from the definition
that the functor $\wt{\Phi}_{P}$ commutes with tensoring by characters of $\G_m$. Hence,
$\Phi_P(F)\simeq\Phi_P(F\{i\})$.
It follows that the map $\phi_P$ 
is $\Z/d$-invariant, which is equivalent to the equality in question.
\ed

\subsection{Orlov's equivalence}\label{Orlov-sec}

Now let $W(x_1,\ldots,x_n)$ be a {\it homogeneous} polynomial of degree $d>0$ (so $\deg(x_i)=1$)
with an isolated singularity.
Let us recall Orlov's construction in \cite[Sec.\ 2.1, 2.2]{Orlov-sing} 
relating the category of $\G_m$-equivariant matrix factorizations of $W$ with the derived category of
coherent sheaves on the smooth projective hypersurface $Y\sub\P^{n-1}$ given by the equation $W=0$.

The construction proceeds in several steps. First, consider the graded algebra 
$$A=k[x_1,\ldots,x_n]/(W).$$
Note that the algebra $A$ is Gorenstein with the Gorenstein parameter $a=n-d$,
i.e., $\Ext^*_A(k,A)$ is concentrated in internal degree $a$ (and in cohomological degree $n-1$).
We denote by $\gr-A$ the category of finitely generated graded $A$-modules, and by
$\gr-A_{\ge i}$ the full subcategory of modules $M$ with $M_j=0$ for $j<i$.

We have Serre's description of coherent sheaves on $Y=\Proj(A)$ as the quotient 
$$\qgr A=\gr-A/\tors-A$$
of $\gr-A$ by the subcategory of torsion modules.
Thus, we have an equivalence
$$D^b(Y)\simeq D^b(\qgr A).$$
Note that under this equivalence the sheaf $\OO_Y(i)$ corresponds to $A(i)\in \qgr A$.

Next, for each $i\in\Z$ we have a fully faithful functor
$$R\om_i: D^b(\qgr A)\to D^b(\gr-A_{\ge i})\sub D^b(\gr-A),$$
which is right adjoint to the natural projection $\pi_i:D^b(\gr-A_{\ge i})\to D^b(\qgr A)$.
The image of $R\om_i$ is denoted by $\DD_i$. The natural projection
$$\pi: D^b(\gr-A)\to D^b(\qgr A)\simeq D^b(Y)$$ 
induces an equivalence $\DD_i\simeq D^b(Y)$.

Recall that a full triangulated subcategory $\TT'$ of a triangulated category $\TT$ is called {\it left (resp., right) admissible}
if the inclusion functor $\TT'\to \TT$ has a left (resp., right) adjoint functor $\TT\to\TT'$. A subcategory is {\it admissible} if 
it is left and right admissible.
A {\it semiorthogonal decomposition}
\begin{equation}\label{semiorthogonal-notation}
\TT=\lan\AA_1,\ldots,\AA_r\ran
\end{equation}
is given by a collection of full triangulated subcategories $(\AA_i)$ such that there exists an increasing filtration
$0=\TT_0\sub\TT_1\sub\ldots\sub\TT_r=\TT$ by left admissible subcategories such that $\AA_i$ is the left orthogonal
of $\TT_{i-1}$ in $\TT_i$. In particular, $\Hom_{\TT}(\AA_j,\AA_i)=0$ for $j>i$ and $(\AA_i)$ generate $\TT$ as a triangulated 
category. Note that if in addition $\TT=D^b(X)$, where $X$ is a smooth projective variety, then each $\AA_i$ is an admissible subcategory in $\TT$ (see \cite{BK}). When one of the subcategories $\AA_i$ is generated by an exceptional object $A_i$
then we write simply $A_i$ instead of $\AA_i$ in the right-hand side of \eqref{semiorthogonal-notation}.

Let $\PP_{\ge i}$ (resp., $\SS_{<i}$) denote the triangulated subcategory of $D^b(\gr-A)$
generated by $A(e)$ with $e\le -i$ (resp., by $k(e)$ with $e>i$).
Then we have semiorthogonal decompositions
$$D^b(\gr-A)=\lan \SS_{<i}, D^b(\gr-A_{\ge i})\ran,$$
$$D^b(\gr-A_{\ge i})=\lan \PP_{\ge i},\TT_i\ran,$$
where $\TT_i$ is equivalent via the natural projection
from $D^b(\gr-A_{\ge i})$ to the graded singularity category of $A$, which in turn is equivalent
to the homotopy category of graded matrix factorizations (by \cite[Thm.\ 3.10]{Orlov-sing}).
Combining these two decompositions we get a semiorthogonal decomposition
\begin{equation}\label{T-i-main}
D^b(\gr-A)=\lan\SS_{<i},\PP_{\ge i},\TT_i \ran.
\end{equation}
On the other hand, there is a semiorthogonal decomposition
\begin{equation}\label{D-i-main}
D^b(\gr-A)=\lan \PP_{\ge i+a},\SS_{<i},\DD_i \ran,
\end{equation}
where $a=n-d$.

For $a=0$ the above decompositions for $i=0$
imply that $\DD_0=\TT_0$, and hence we get an equivalence relating matrix factorizations and
coherent sheaves on $Y$:
$$\DMF_{\G_m}(W)\simeq \TT_0=\DD_0\simeq D^b(Y).$$

More generally, for $a\ge 0$ we can use the semiorthogonal decomposition
$$\PP_{\ge i}=\lan \PP_{\ge i+a}, A(-i-a+1),\ldots, A(-i-1), A(-i)\ran$$
to refine the decomposition \eqref{T-i-main} to
$$D^b(\gr-A)=\lan \SS_{<i}, \PP_{\ge i+a}, A(-i-a+1),\ldots, A(-i-1), A(-i), \TT_i\ran.$$
Comparing this with \eqref{D-i-main} we get a semiorthogonal
decomposition
\begin{equation}\label{D-i-T-i-eq}
\DD_i=\lan A(-i-a+1),\ldots, A(-i-1), A(-i), \TT_i\ran.
\end{equation}
Thus, in this case the category $D^b(Y)$ contains $\DMF_{\G_m}(W)$ as an admissible subcategory.
Let
\begin{equation}\label{rho-1-eq}
\rho:D^b(Y)\simeq \DD_0\to \TT_0\simeq\DMF_{\G_m}(W)
\end{equation}
be the right adjoint functor to the embedding $\TT_0\sub \DD_0$.

In the case $a\le 0$, using the decomposition
$$\SS_{<-a}=\lan \SS_{<0}, k, k(-1), \ldots, k(a+1)\ran,$$
we refine the decomposition \eqref{D-i-main} for $i=-a$ to
$$D^b(\gr-A)=\lan \PP_{\ge 0}, \SS_{<0}, k, k(-1), \ldots, k(a+1), \DD_{-a} \ran.$$
Comparing this with \eqref{T-i-main} for $i=0$ we get a semiorthogonal
decomposition
\begin{equation}\label{MF-decomposition-eq}
\TT_0=\lan k, k(-1), \ldots, k(a+1), \DD_{-a} \ran.
\end{equation}
Thus, in this case $\DMF_{\G_m}(W)$ contains $D^b(Y)$ as a full subcategory.
Let
\begin{equation}\label{rho-2-eq}
\rho:D^b(Y)\simeq \DD_{-a}\to \TT_0\simeq\DMF_{\G_m}(W)
\end{equation}
be the corresponding fully faithful functor.

Caldararu and Tu showed that the above constructions can be performed at the dg-level (see \cite[Sec.\ 5]{CT}).
In particular, the functors \eqref{rho-1-eq} and \eqref{rho-2-eq} lift to the dg-level.

For any graded module $M$ over $A$ we denote by $M^{\st}$ the object of $\DMF_{\G_m}(W)$ corresponding
to $M$ viewed as an object of the graded singularity category.

\begin{prop}\label{Orlov-prop} 
Let us consider the composition
$$\Phi:D^b(\P^{n-1})\to D^b(Y)\rTo{\rho} \TT_0\simeq \DMF_{\G_m}(W),$$
where the first arrow is the pull-back with respect to the embedding $Y\sub \P^{n-1}$, and
$\rho$ is given by \eqref{rho-1-eq} for $a\ge 0$ and by \eqref{rho-2-eq} for $a\le 0$.
Then the image of $\Phi$ is contained in the triangulated subcategory of
$\DMF_{\G_m}(W)$, generated by $(k(-e)^{\st})_{e\in\Z}$.
\end{prop}

\Pf . Assume first that $a\le 0$.
For $e<0$ let us consider the truncated module $A(-e)_{\ge -a}$, so that we have an exact
sequence
\begin{equation}\label{trunc-ex-seq}
0\to A(-e)_{\ge -a}\to A(-e)\to N_e\to 0,
\end{equation}
where $N_e\in\lan k(-e), k(-e-1),\ldots, k(a+1)\ran$.
Then we claim that $A(-e)_{\ge -a}\in \DD_{-a}$. Indeed, we have
$A(-e)_{\ge -a}\in D^b(\gr-A_{\ge -a})$, so by the decomposition \eqref{D-i-main}, it is enough to check that
$$\Ext^*_{\gr-A}(A(-e)_{\ge -a},\PP_{\ge 0})=0.$$
But this follows from the exact sequence \eqref{trunc-ex-seq}, since
$\Ext^*_{\gr-A}(k(-j),\PP_{\ge 0})=0$ for $j<-a$.
Now the same exact sequence shows that on the one hand,
$$\pi(A(-e)_{\ge -a})=\OO(-e),$$
while on the other hand, the image of $A(-e)_{\ge -a}$ in the graded singularity category is
the same as that of $N_e[-1]$. Hence, we deduce that 
$$\Phi(\OO(-e))\in \lan k(-e)^{\st}, k(-e-1)^{\st},\ldots, k(a+1)^{\st}\ran\sub\DMF_{\G_m}(W).$$
Since $D^b(\P^{n-1})$ is generated by the sheaves $\OO(-e)$ with $e<0$, the assertion follows in this case.

Now assume that $a\ge 0$. Then for $e<0$ we consider the truncated module
$A(-e)_{\ge 0}$ that fits into an exact sequence
$$0\to A(-e)_{\ge 0}\to A(-e)\to M_e\to 0$$
where $M_e\in\lan k(-e), \ldots, k(1)\ran$.
Then as above we deduce that $A(-e)_{\ge 0}\in \DD_0$ and 
$$\pi(A(-e)_{\ge 0})=\OO(-e).$$ 
On the other hand, by the semiorthogonal decomposition \eqref{D-i-T-i-eq} for $i=0$, we
have an exact triangle
$$\rho(A(-e)_{\ge 0})\to A(-e)_{\ge 0}\to Q\to \ldots$$
with $Q\in \lan A(-a+1),\ldots, A(-1), A\ran$. It follows that the image of $\rho(A(-e)_{\ge 0})$
in the graded singularity category is the same as that of $M_e$, so we deduce that
$$\Phi(\OO(-e))\in \lan k(-e)^{\st}, \ldots, k(1)^{\st}\ran\sub\DMF_{\G_m}(W).$$
\ed

\section{Rational structure on the Hochschild homology}
\label{rat-lattice-sec}

\subsection{Rational structure on the Hochschild homology of admissible subcategories in the derived categories of sheaves}
\label{rational-lattice-adm-sub-sec}

Let $X$ be a smooth projective variety over $\C$. Recall that the Hochschild homology $HH_*(X)$
can be identified with $H^*(X,\C)=\bigoplus_{p,q} H^{p,q}(X)$ (with $H^{p,q}(X)\sub HH_{p-q}(X)$). 
We can use this identification to define a rational lattice
in $HH_*(X)$. To get better compatibility with the Chern characters and Fourier-Mukai type functors we insert some standard factors. 
Namely, let us consider an automorphism 
$$J:H^*(X,\C)\to H^*(X,\C): c\mapsto (2\pi i)^p c \ \ \text{ for } c\in H^{p,q}(X),$$
and set\footnote{Choosing  instead of $J$ the automorphism $c\mapsto (2\pi i)^q c$, $c\in H^{p,q}(X)$, would work as well. 
The two choices differ by the grading operator with respect to the Hochschild degree.} 
$$HH_*(X)_\Q:=J(H^*(X,\Q))\sub H^*(X,\C)\simeq HH_*(X).$$

\begin{prop}\label{rational-prop} 
(i) For any $E\in D^b(X)$ we have $\ch(E)\in HH_*(X)_\Q$.

\noindent
(ii) For any $c,c'\in HH_*(X)_\Q$ one has $(c,c')_X\in \Q$, where $(\cdot,\cdot)_X$ is the pairing \eqref{X-pairing}.

\noindent
(iii) For any Fourier-Mukai type functor $\Phi:D^b(X)\to D^b(Y)$, where 
$Y$ is smooth projective, the induced map
$\Phi_*:HH_*(X)\to HH_*(Y)$ sends $HH_*(X)_\Q$ to $HH_*(Y)_\Q$.
\end{prop}

\Pf . (i) This follows immediately from Proposition \ref{Chern-prop}.

\noindent
(ii) We have 
$$(c,c')_X=\frac{1}{(2\pi i)^n}\int^{top}_{[X]}c\cdot c'\cdot\Td_X,$$
so the assertion follows from the fact that $c\cdot c'\cdot\Td_X\in J(H^*(X,\Q))$.

\noindent
(iii) This follows from (i) and (ii) and from Lemma \ref{FM-HH-lem}.
\ed

We will use the real structure on $HH_*(X)=H^*(X,\C)$ associated with the rational lattice $HH_*(X)_\Q=J(H^*(X,\Q))$.
Let us denote by $\tau$ the corresponding complex conjugation map on $H^*(X,\C)$, so that
$\tau(Jx)=J\ov{x}$ (where $x\mapsto\ov{x}$ is the usual complex conjugation on $H^*(X,\C)$). 
It is easy to check that 
\begin{equation}\label{tau-formula}
\tau(c)=(-1)^p(2\pi i)^{q-p}\ov{c} \ \ \text{ for } c\in H^{p,q}(X).
\end{equation}

Let $\TT\sub D^b(X)$ be an admissible subcategory. We can enhance $\TT$ to a dg-category using the dg-enhancement
of $D^b(X)$ (see \cite[Sec.\ 4]{Kuz}). 
By functoriality of the Hochschild homology, applied to the inclusion functor and to its left adjoint
$D^b(X)\to\TT$, the Hochschild homology $HH_*(\TT)$ gets identified with a direct summand of $HH_*(X)$.

\begin{prop}\label{adm-rat-lattice-prop} 
The subgroup $HH_*(\TT)\cap HH_*(X)_\Q$ is a rational lattice in $HH_*(\TT)$.
If $\TT\hra D^b(X')$ is a different embedding of $\TT$ as admissible subcategory,
where $X'$ is a smooth projective variety, then 
$$HH_*(\TT)\cap HH_*(X)_\Q=HH_*(\TT)\cap HH_*(X')_\Q,$$
so that the lattice $HH_*(\TT)_\Q:=HH_*(\TT)\cap HH_*(X)_\Q$ depends only on $\TT$.
If $\TT\to \TT'$ is a dg-functor between two categories like this then the induced map
$HH_*(\TT)\to HH_*(\TT')$ is compatible with these rational lattices. 
\end{prop}

\Pf . The projector functor $\Pi: D^b(X)\to\TT\sub D^b(X)$ is given by some kernel (see \cite[Thm.\ 7.1]{Kuz-bc}).
Hence, by Proposition \ref{rational-prop}, the induced projector $\Pi_*$ of $HH_*(X)$ with the image $HH_*(\TT)$
sends $HH_*(X)_\Q$ to $HH_*(\TT)\cap HH_*(X)_\Q$, which implies that 
the latter subgroup is a rational lattice in $HH_*(\TT)$.
Suppose we have a functor $\Phi:\TT\to\TT'$, where $\TT'\sub D^b(X')$. Then by Proposition \ref{rational-prop}(iii),
the composed map 
$$HH_*(X)\rTo{\Pi_*} HH_*(\TT)\rTo{\Phi_*} HH_*(\TT')\to HH_*(X')$$
is compatible with rational lattices, hence, $\Phi_*$ sends
$HH_*(\TT)\cap HH_*(X)_\Q=\Pi_*(HH_*(X)_\Q)$ to $HH_*(\TT')\cap HH_*(X')_\Q$.
In the case $\Phi=\Id$ this also proves the independence of the lattice on the embedding $\TT\hra D^b(X)$.
\ed

\subsection{Hodge-Riemann relations for matrix factorizations: homogeneous case}\label{LG-rat-lattice-sec}

We use the rational lattices considered above to define a rational lattice in the Hochschild homology
of the category $\MF_{\G_m}(W)$, where $W$ is a homogeneous polynomial with an isolated singularity.


\begin{lem}\label{MF-admissible-lem} 
Let $W$ be a homogeneous polynomial with isolated singularity.
Then there exists a smooth projective variety $X$ such that
$\DMF_{\G_m}(W)$ is an admissible subcategory in $D^b(X)$,
in a way compatible with the dg-enhancements.
\end{lem}

\Pf . Let $a=n-d$ be the Gorenstein parameter. If $a\ge 0$ then the semiorthogonal decomposition
\eqref{D-i-T-i-eq} shows that we can take $X$ to be the hypersurface $Y\sub\P^{n-1}$
with the equation $W=0$. In the case $a\le 0$ we have the semiorthogonal decomposition \eqref{MF-decomposition-eq}
of $\DMF_{\G_m}(W)$, with $D^b(Y)$ as one of the pieces, where each of the remaining pieces is generated by an
exceptional object. Hence, the desired $X$ can be constructed using \cite[Thm.\ 4.15]{Orlov-adm}. 
\ed

Combining this Lemma with Proposition \ref{adm-rat-lattice-prop} we 
equip the Hochschild homology $\HH(W)=HH_*(\MF_{\G_m}(W))$ with a (uniquely defined) rational lattice,
such that the maps on Hochschild homology induced by Fourier-Mukai transforms involving $\MF_{\G_m}(W)$ 
are compatible with this lattice.
We denote by $x\mapsto \ov{x}$ the complex conjugation associated with the corresponding real structure on $\HH(W)$.

\begin{lem}\label{rat-subspace-lem} 
The subspace $H(W)^{\mu_d}\sub \HH(W)$ coming from the decomposition \eqref{HH-W-decomp}
is compatible with the rational lattice in $\HH(W)$. The projector $\Pi: \HH(W)\to H(W)^{\mu_d}$ is compatible
with the rational lattices.
\end{lem}

\Pf . Indeed, the operators of the $\Z/d$-action on $\HH(W)$ are induced by the twist functors $E\mapsto E\{m\}$, hence
they are compatible with the rational lattice in $\HH(W)$. It remains to observe that 
$H(W)^{\mu_d}$ is the subspace of $\Z/d$-invariants in $\HH(W)$, and
$\Pi$ is precisely the standard projector onto it (see \cite[Thm.\ 2.6.1(ii)]{PV-CohFT}). 
\ed

Let us denote by $\lan\cdot,\cdot\ran_W$ the canonical pairing \eqref{pair-eq} on $\HH(W)$.

\begin{lem}\label{W-orth-lem}
For any $x\in H(W)^{\mu_d}\sub\HH(W)$ and any $m\in\Z$ one has
$$\lan x,\ch(k(m)^{\st})\ran_W=0.$$
\end{lem}

\Pf . This follows from the explicit calculation of the pairing $\lan\cdot,\cdot\ran_W$ and of $\ch(k(m)^{\st})$ in \cite{PV-HRR}.
Note that the $\Z/2$-folding of the dg-category $\MF_{\G_m,d}(W)$ is naturally isomorphic to the $\Z/2$-dg-category
$\MF_{\mu_d}(W)$ of $\mu_d$-equivariant matrix factorizations (see \cite[Sec.\ 2.1]{PV-CohFT} and \cite[Sec.\ 4.4]{PV-HRR}). The computations
in \cite{PV-HRR} were done in the context of $\Z/2$-dg-categories but we can use them due to the compatibility of
all the Hochschild homology manipulations with the $\Z/2$-folding (see remarks at the end of Sec.\ \ref{dg-sec}).
By \cite[Prop.\ 4.3.4]{PV-CohFT}, we see that $\ch(k^{\st})$ has a trivial component in $H(W)^{\mu_d}$, i.e., $\Pi(\ch(k^{\st}))=0$.
It follows that $\Pi(\ch(k(m)^{\st}))=0$ for any integer $m$ (since $\Pi$ is the projector onto the invariants of the $\Z/d$-action). 
Now the assertion follows from the fact that the summand $H(W)^{\mu_d}$ is orthogonal to other summands in the decomposition \eqref{HH-W-decomp} with respect to $\lan\cdot,\cdot\ran_W$, as the explicit formula of
\cite[Thm.\ 4.2.1]{PV-HRR} shows.
\ed

\begin{prop}\label{HR-prop} 
Suppose we have a class $x\in H(W)^{\mu_d}\cap \HH_j(W)$ for some $j\in\Z$. If $\lan x,\ov{x}\ran_W=0$ then
$x=0$.
\end{prop}

\Pf . Let $a$ be the Gorenstein parameter. Assume first that $a\ge 0$. Then we have a fully faithful functor
$\la:\DMF_{\G_m}(W)\to D^b(Y)$ and the right adjoint functor $\rho$ (see \eqref{rho-1-eq}), such that $\rho\la=\Id$.
By Proposition \ref{Orlov-prop}, the map
$$\rho_*:H^*(Y,\C)\to \HH(W)$$
sends classes restricted from $\P^{n-1}$ to the span of the Chern characters $\ch(k(m)^{\st})$, $m\in\Z$.
But the latter classes are orthogonal to $H(W)^{\mu_d}$ with respect to $\lan\cdot,\cdot\ran_W$ 
by Lemma \ref{W-orth-lem}. Hence,
by adjointness of $\la_*$ and $\rho_*$ (see Lemma \ref{adjoint-HH-lem}),
$\la_*(H(W)^{\mu_d})$ is left orthogonal to the image of $H^*(\P^{n-1},\C)\to H^*(Y,\C)$ with respect to the pairing
$\lan\cdot,\cdot\ran_{D^b(Y)}$. Since $\Td_Y$ is a class restricted from $\P^{n-1}$, using \eqref{Ram-formula} we get that
$$\int_Y \kappa(\la_*(H(W)^{\mu_d}))\cdot c=0$$
for any $c$ restricted from $\P^{n-1}$.
Hence, for any $x\in H(W)^{\mu_d}$, $\la_*(x)$ is a primitive class. In particular, by the Lefschetz hyperplane theorem,
$\la_*(x)\in H^{n-2}(Y,\C)$.
If in addition $x\in \HH_j(W)$ then $\la_*(x)\in H^{p,q}(Y)$ for the unique $p$, $q$ such that
$p+q=n-2$ and $p-q=j$. Since $\la_*$ is compatible with rational lattices, we have
$$\la_*(\ov{x})=\tau(\la_*(x))=(-1)^p(2\pi i)^{q-p}\ov{\la_*(x)}$$
(see \eqref{tau-formula}). Hence,
$$(-1)^p(2\pi i)^{q-p}\lan \la_*(x),\ov{\la_*(x)}\ran_{D^b(Y)}=\lan \la_*(x),\la_*(\ov{x})\ran_{D^b(Y)}
=\lan x,\rho_*\la_*(\ov{x})\ran_W=\lan x,\ov{x}\ran_W.$$
Thus, the vanishing of $\lan x,\ov{x}\ran_W$ implies the vanishing of $\lan \la_*(x),\ov{\la_*(x)}\ran_{D^b(Y)}$.
Since $\la_*(x)$ is primitive, by Lemma \ref{HR-lem}, this implies that $\la_*(x)=0$, and so $x=0$.

Now assume that $a\le 0$.
Then we have a fully faithful functor 
$$\rho:D^b(Y)\to\DMF_{\G_m}(W)$$
(see \eqref{rho-2-eq}), such that the image is the left orthogonal to $k^{\st},\ldots,k(a+1)^{\st}$. Considering
the corresponding decomposition of $\HH(W)$ (see \cite[Thm.\ 7.3]{Kuz}) we deduce
that the image of $\rho_*$ can be identified with the left orthogonal to $\ch(k^{\st}),\ldots,\ch(k(a+1)^{\st})$
with respect to $\lan\cdot,\cdot\ran_W$. Since $H(W)^{\mu_d}$ lies in this left orthogonal, we deduce that
any $x\in H(W)^{\mu_d}\cap \HH_j(W)$ has form $x=\rho_*(y)$ for some $y\in HH_j(Y)$. 
Furthermore, by Proposition \ref{Orlov-prop} and 
Lemma \ref{W-orth-lem}, we have
$$\lan y,c\ran_{D^b(Y)}=\lan x,\rho_*(c)\ran_W=0$$
for any $c$ restricted from $\P^{n-1}$. Thus, we deduce that $y$ is a primitive class. Now
we can finish the proof as before, using Lemma \ref{HR-lem}
and the fact that $\rho_*$ is compatible with the rational lattices.
\ed

\begin{rem}\label{CY-rem} 
It is easy to see that the Chern characters $\ch(k(m)^{\st})$, $m\in\Z$, span
the orthogonal complement to $H(W)^{\mu_d}$ in $\HH(W)$.
In the Calabi-Yau case, $d=n$, Proposition \ref{Orlov-prop} implies that the subspace
$H(W)^{\mu_d}\sub \HH(W)$ corresponds to the primitive part of the middle cohomology of the projective hypersurface $Y$
under the isomorphism
$$\HH(W)\simeq H^*(Y,\C)$$
induced by the Orlov's equivalence. Note that the images of $k(m)^{\st}$ in $D^b(Y)$ under this equivalence
are calculated explicitly in \cite[Prop.\ 4.11]{CIR}.
\end{rem}

\subsection{Hodge-Riemann relations for matrix factorizations: quasihomogeneous case}\label{LG-rat-lat-qh-sec}

Now let $W(x_1,\ldots,x_n)$ be a quasihomogeneous polynomial with an isolated singularity, such that the corresponding
homogeneous polynomial $\wt{W}(y_1,\ldots,y_n)=W(y_1^{d_1},\ldots,y_n^{d_n})$ still has an isolated singularity.
Let us consider the corresponding finite flat $\G_m$-equivariant morphism between affine spaces
$$\varphi:\A^n\to \A^n:(y_1,\ldots,y_n)\mapsto (y_1^{d_1},\ldots, y_n^{d_n}),$$
such that $\varphi^*W=\wt{W}$. We have the corresponding functors
$$\varphi^*:\MF_{\G_m}(W)\to \MF_{\G_m}(\wt{W}), \ \ \varphi_*:\MF_{\G_m}(\wt{W})\to \MF_{\G_m}(W),$$
such that $\varphi_*\varphi^*(E)\simeq \varphi_*\OO\ot E$, where $\varphi_*\OO$ corresponds to a free $\C[x_1,\ldots,x_n]$-module
(with generators of various degrees).

\begin{lem}\label{quasihom-lattice-lem}
The induced maps on the Hochschild homology
$$\varphi^*:\HH(W)\to \HH(\wt{W}), \ \  \varphi_*:\HH(\wt{W})\to \HH(W)$$ 
are compatible with the decompositions \eqref{HH-W-decomp}, and the composition
$$H(W)^{\mu_d}\rTo{\varphi^*} H(\wt{W})^{\mu_d} \rTo{\varphi_*} H(W)^{\mu_d}$$
is the multiplication by $\deg(\varphi)=d_1\ldots d_n$. 
\end{lem}

\Pf . The first assertion immediately follows from the fact that $\varphi^*$ and $\varphi_*$ commute with the twist functors
$E\mapsto E\{m\}$. Next, since $\varphi_*\OO$ is free of rank $\deg(\varphi)$,
 we see that $\varphi_*\varphi^*(E)$ is a direct sum of $\deg(\varphi)$ twists $E\{m\}$. It remains to use the fact that these twists
 act trivially on the summand $H(W)^{\mu_d}$.
\ed

The above Lemma shows that the map
$$\varphi^*:H(W)^{\mu_d}\to H(\wt{W})^{\mu_d}$$
is injective and its image coincides with the image of the map $\varphi^*\varphi_*:H(\wt{W})^{\mu_d}\to H(\wt{W})^{\mu_d}$.
Recall that by Lemma \ref{rat-subspace-lem}, the subspace $H(\wt{W})^{\mu_d}\sub \HH(\wt{W})$ is compatible
with the rational lattice in $\HH(\wt{W})$. Since $\varphi^*\varphi_*$ is induced by a dg-endofunctor of $\MF_{\G_m}(\wt{W})$,
it follows that the subspace
$$\varphi^*(H(W)^{\mu_d})=\varphi^*\varphi_*(H(\wt{W})^{\mu_d})\sub H(\wt{W})^{\mu_d}$$
inherits a rational lattice, so we get a rational lattice on $H(W)^{\mu_d}$.

\begin{lem}\label{quasihom-lattice-fun-lem} 
For any dg-functors $F: \DMF_{\G_m}(W)\to \TT$, $G:\TT\to \DMF_{\G_m}(W)$, where $\TT$ is an admissible subcategory
in $D^b(X)$ for some smooth and projective $X$, the induced maps
$$H(W)^{\mu_d}\to \HH(W)\rTo{F_*} HH_*(\TT), \ \ HH_*(\TT)\rTo{G_*} \HH(W)\rTo{\Pi} H(W)^{\mu_d}$$
are compatible with rational structures.
\end{lem}

\Pf . By Lemma \ref{quasihom-lattice-lem}, to prove the assertion about $F_*$ we can replace it by
$F_*\varphi_*\varphi^*$. Since $F_*\varphi_*$ is induced by a dg-functor $\DMF_{\G_m}(\wt{W})\to\TT$, it is 
compatible with the rational lattices. But the restriction of $\varphi^*$ to $H(W)^{\mu_d}$ preserves rational lattices by
the definition, so the assertion follows.

To check the assertion about $\Pi G_*$, it is enough to prove it for the composition 
$$\varphi^*\Pi G_*=\Pi \varphi^*G_*: HH_*(\TT)\to H(\wt{W})^{\mu_d}.$$
But $\varphi^*G_*$ is induced by a dg-functor $\TT\to \DMF_{\G_m}(\wt{W})$, so it compatible with the rational lattices.
It remains to use the fact that $\Pi:\HH(\wt{W})\to H(\wt{W})^{\mu_d}$ is also compatible with the rational lattices
(see Lemma \ref{rat-subspace-lem}). 
\ed

\begin{rem} Using the connection between $\DMF_{\G_m}(W)$ and the derived category of the corresponding stacky
weighted projective hypersurface $\XX$ 
from \cite{Orlov-sing} and the recent paper \cite{BLS},
one can equip the Hochschild homology 
$\HH(W)$ with a rational structure for any quasihomogeneous polynomial $W$ with an isolated singularity. 
Namely, the main result of \cite{BLS} implies that 
$D^b(\XX)$ can be realized as an admissible subcategory in $D^b(Y)$ for $Y$ smooth and projective variety, so using the same
approach as in Lemma \ref{MF-admissible-lem} we can realize $\DMF_{\G_m}(W)$ as such a subcategory.
\end{rem}

As before, we denote by $x\mapsto \ov{x}$ the conjugation associated with the real structure on $H(W)^{\mu_d}$.

\begin{prop}\label{quasihom-HR-prop} Let $W$ be a quasihomogeneous polynomial with an isolated singularity
such that $\wt{W}$ still has an isolated singularity.
Given a class $x\in H(W)^{\mu_d}\cap \HH_j(W)$ for some $j\in\Z$, if $\lan x,\ov{x}\ran_W=0$ then
$x=0$.
\end{prop}

\Pf . The map $\varphi^*:H(W)^{\mu_d}\to H(\wt{W})^{\mu_d}$ is rational, so using adjointness of $(\varphi^*,\varphi_*)$ and 
Lemma \ref{adjoint-HH-lem} we get
$$\lan \varphi^*x, \ov{\varphi^*(x)}\ran_{\wt{W}}=\lan \varphi^*x,\varphi^*(\ov{x})\ran_{\wt{W}}=\lan x, \varphi_*\varphi^*(\ov{x})\ran_W=\deg(\varphi)\lan x, \ov{x}\ran_W=0.$$
Also we have $\varphi^*x\in H(\wt{W})^{\mu_d}\cap \HH_j(\wt{W})$.
Hence, by Proposition \ref{HR-prop}, we obtain $\varphi^*x=0$, and so by Lemma \ref{quasihom-lattice-lem}, $x=0$.
\ed

\begin{rem} The idea of using the relation between matrix factorizations of $W$ and $\wt{W}$ to deduce results about
$\HH(W)$ was inspired by a similar method in \cite[Sec.\ 6]{MPSW}.
\end{rem}

\section{Homogeneity}\label{homog-sec}

\subsection{Calculations with Koszul matrix factorizations}\label{Koszul-mf-calc-sec}

Let $A$ and $B$ be $\G_m$-vector bundles over a scheme $S$ (where $\G_m$ acts trivially on $S$),
$\a\in H^0(\tot(A),p^*B^\vee\{d\})$ and $\b\in H^0(\tot(A),p^*B)$
be $\G_m$-invariant sections, where $p:\tot(A)\to S$ is the projection.
Assume that $\a$ and $\b$ are orthogonal and have common zeros only on the zero section in $\tot(A)$.
Let 
$$E=\{\a,\b\}\in\MF_{\G_m,d}(\tot(A),0)$$ 
be the corresponding Koszul matrix factorization of $0$ on $\tot(A)$ (see Section \ref{mf-sec}).
The following homogeneity property follows from the results of \cite{PV-Wc} and \cite{Chiodo}
(it also appears implicitly in the proof of \cite[Prop.\ 5.6.1]{PV-CohFT}).

\begin{prop}\label{homog-class-prop} 
In the above situation we have
$$\Td(B)\Td(A)^{-1}\ch(\ov{\com}(p_*E))\in H^{2(\rk B-\rk A)}(S).$$
\end{prop}

\Pf . Note that 
$$\ch(\ov{\com}(p_*E))=\ch(p_*E):=\ch(H^{even}(p_*E))-\ch(H^{odd}(p_*E)),$$ 
where on the right we view $p_*E$ as a $\Z/2$-graded complex.
It is enough to show that in the Chow group $A^*(S)_\Q$ one has
$$\Td(B)\Td(A)^{-1}\ch(p_*E)\in A^{\rk B-\rk A}(S)_\Q.$$
By \cite[Lemma 5.3.8]{Chiodo}, one has
$$\Td(A)^{-1}\ch(p_*E)=\ch_S^{\tot{A}}(E)\cdot [p],$$
where $\ch_S^{\tot{A}}(E)\in A^*(S\to \tot(A))$
is the localized Chern character of the $\Z/2$-graded complex $E$ in the relative Chow group 
for the zero section embedding $S\to \tot(A)$ (see \cite[Sec.\ 2.2]{PV-Wc}),
and $[p]\in A^{-rk A}(\tot{A}\to S)$ is the orientation class of $p$. Now by
\cite[Thm.\ 3.2]{PV-Wc}, the class
$$\Td(B)\ch_S^{\tot{A}}(E)\in A^*(S\to \tot{A})$$
is concentrated in degree $\rk B$. To get the class we need, one has to multiply the above class with the orientation class
$[p]$ that lives in degree $-rk A$, hence the result.
\ed

Now let us consider the setup of Sec.\ \ref{purity-sec}. Note that the functor $\Phi:\DMF_{\G_m}(W)\to D^b(X)$ 
in \eqref{main-Phi-eq} is of the form $\Phi=\Phi_P$ (see \eqref{Phi-P-eq}) for the diagram
\begin{diagram}
&&\tot(A)\\
&\ldTo{Z}&&\rdTo{p}&\\
\A^n&&&& X
\end{diagram}
and $P=\{\a,\b\}$. Note also that in this case $\DMF_{\G_m}(W)=\ov{\DMF}_{\G_m}(W)$
is the usual homotopy category of matrix factorizations.

We start by computing the left adjoint functor $\Psi:D^b(X)\to \DMF_{\G_m}(W)$ to $\Phi$.

\begin{lem} In the notation of \eqref{Phi-P-eq}, 
one has $\Psi=\Psi_{P'}$, where 
$$P'=\{ \a',\b'\}\ot p^*({\det}^{-1}(A)\ot\om_X)[N],$$ 
with
$\a'=\b\in p^*B$, $\b'=-\a\in p^*B^\vee\{d\}$, $N=\dim X+\rk A-n$. 
Thus, the Koszul matrix factorization $\{\a',\b'\}$ is of the same type as 
$\{\a,\b\}$ but with $B$ being replaced by $B^\vee\{d\}$.
\end{lem}

\Pf . By Proposition \ref{mf-adjoint-prop}, we have $\Psi=\Psi_{P'}$ with
$$P'=P^\vee\ot \om_Z[N].$$
Recall that the dual matrix factorization $P^\vee$ has the even part $P_0^\vee$, the odd part $P_1^\vee\{-d\}$, and
the differential $\de_{P^\vee}$ determined by the rule
$$\lan \de_{P^\vee}(\xi),x\ran=(-1)^{\deg(\xi)}\lan \xi,\de_P(x)\ran.$$
Thus, disregarding the $\G_m$-action we can identify 
$P^\vee$ with ${\bigwedge}^*(p^*B)$ using the isomorphism
$${\bigwedge}^* (p^*B)\rTo{\sim} ({\bigwedge}^*(p^*B^\vee))^\vee:
b_1\we\ldots\we b_p\mapsto [w\mapsto (\iota_{b_1}\circ\ldots\circ\iota_{b_p})(w)_0],$$
where the subscript $0$ means taking the component in $\bigwedge^0$. Under this identification,
with the above sign convention, the operator dual to $\a\we$ is $-\iota_{\a}$, and the operator dual to 
$\iota_\b$ is $\b\we$. The $\G_m$-weights work out so that we have an identification
$$P^\vee\simeq\{\a',\b'\}.$$
Taking into account the isomorphism 
$$\om_Z\simeq \om_{\tot A}\simeq \om_p\ot p^*\om_X\simeq p^*({\det}^{-1}(A)\ot\om_X),$$ 
we get the assertion.
\ed

\begin{lem}\label{Koszul-comp-lem}
In the above situation the composed functor $\Phi\circ\Psi: D^b(X)\to D^b(X)$
is associated with a kernel $K\in D^b(X\times X)$ such that $[K]=d\cdot[K']$ in $K_0(X\times X)$, where
$$K'=\ov{\com}\left(p^{(2)}_*\{\a^{(2)},\b^{(2)}\}\right)\ot p_1^*({\det}^{-1}(A)\ot\om_X)[N]$$ 
for a Koszul
matrix factorization of zero $\{\a^{(2)},\b^{(2)}\}$ on the vector bundle 
$$p^{(2)}:\tot(A)\times_{\A^n}\tot(A)\to X\times X,$$
associated with the sections $\a^{(2)}=(\a,\a')\in p_1^*B^\vee\{d\}\oplus p_2^*B$, 
$\b^{(2)}=(\b,\b')\in p_1^*B\oplus p_2^*B^\vee\{d\}$.
Furthermore, $\a^{(2)}$ and $\b^{(2)}$ have common zeros only on the zero section of the bundle $p^{(2)}$.
\end{lem}

\Pf . The computation of the composition follows from Lemma \ref{gen-comp-lem}. The fact about the common zeros
of $\a^{(2)}$ and $\b^{(2)}$ follows from the similar fact about $(\a,\b)$ and $(\a',\b')$.
\ed

Note that the bundle $\tot(A)\times_{\A^n}\tot(A)$ over $X\times X$ has rank $2\rk A-n$.
Thus, by Proposition \ref{homog-class-prop}, we have
$$\Td(p_1^*B^\vee)\Td(p_2^*B)\Td(p_1^*A)^{-1}\Td(p_2^*A)^{-1}\ch\bigl(K\ot p_1^*(\det(A)\ot\om_X^{-1})\bigr)\in 
H^{2(2\rk B-2\rk A+n)}(X\times X).$$
Using the formula
$$\Td(E^\vee)=\Td(E)\cdot \ch(\det(E))^{-1},$$
we can rewrite the above class as
\begin{equation}\label{main-homog-eq}
p_1^*\left(\ch(\om_X)^{-1}\Td(B^\vee)\Td(A^\vee)^{-1}\right)\cdot p_2^*\left(\Td(B)\Td(A)^{-1}\right)\ch(K)\in 
H^{2D}(X\times X),
\end{equation}
where $D$ is given by \eqref{D-eq}.


\subsection{The proof of Theorem \ref{main-thm}}


We start by reformulating the statement using the canonical bilinear forms on Hochschild homology.
Let $\phi=\Phi_*:\HH(W)=HH_*(\MF_{\G_m}(W))\to H^*(X)$ be the map induced by $\Phi=\Phi_P$ on Hochschild homology.
Let us set
$$\a=\Td(B)\Td(A)^{-1},$$ 

By the nondegeneracy of the Poincar\'e pairing, and by \eqref{Ram-formula},
the left orthogonal to $H^i(X)\sub H^*(X)$ with respect to $\lan\cdot,\cdot\ran_{D^b(X)}$
is the subspace
$$^{\perp}H^i(X)=\bigoplus_{j\neq 2\dim X-i} \kappa(\Td_X)^{-1}\cdot H^j(X)\sub H^*(X),$$
where $\kappa$ is given by \eqref{iota-eq}.
Since $\lan x,\a^{-1}\cdot y\ran_{D^b(X)}=\lan \kappa(\a^{-1})\cdot x, y\ran_{D^b(X)}$, we deduce that
\begin{equation}\label{left-perp-eq}
^{\perp}\bigl(\a^{-1}\cdot H^i(X)\bigr)=\bigoplus_{j\neq 2\dim X-i} \kappa(\a\cdot\Td_X)^{-1}\cdot H^j(X)\sub H^*(X)
\end{equation}
(note that $\kappa(\a)$ and $\kappa(\Td_X)$ live in even degrees, so they commute with any cohomology class).

By the nondegeneracy of $\lan\cdot,\cdot\ran_{D^b(X)}$, 
to show that $\phi(x)\in \a^{-1}\cdot H^D(X)$ for all $x\in \HH(W)$, it is enough to prove that 
$$\psi\bigl(^{\perp}\bigl(\a^{-1}\cdot H^D(X)\bigr)\bigr)=0,$$ 
where $\psi$ is the left adjoint operator to $\phi$.
Note that by Lemma \ref{adjoint-HH-lem}, we have $\psi=\Psi_*$. 
Thus, taking into account \eqref{left-perp-eq}, we should check that
for each $j\neq 2\dim X-D$ one has
$$\psi \left(\kappa(\a\cdot\Td_X^{-1})\cdot H^j(X)\right)=0.$$

We are going to use the rational lattices on the relevant Hochschild homology introduced in Sec.\ \ref{rat-lattice-sec}.
We denote by $x\mapsto \ov{x}$ (resp., $\tau$) the corresponding operation of complex conjugation on $\HH(W)$
(resp., $H^*(X,\C)$).

Let $y\in H^{p,q}(X)\sub H^j(X)$, and set $y'=\kappa(\a\Td_X^{-1})\cdot y$.
Note that $y'$, viewed as an element of Hochschild homology, lives in the single degree $p-q$.
Hence, the same is true about $z=\psi(y')$.
Also, since by Lemma \ref{phi-inv-lem}, $\phi\Pi=\phi$, it follows that $\Pi\psi=\psi$, so $z\in H(W)^{\mu_d}$.
By Proposition \ref{quasihom-HR-prop}, $z=0$ if and only if $\lan\ov{z},z\ran_W=0$.
Thus, it is enough to prove that
$$\lan\ov{\psi (y')},\psi (y')\ran_W=0,$$
provided $j\neq 2d-D$.

We will use the following fact about the involution $\kappa$. 
For any vector bundle $V$ on $X$ one has
\begin{equation}\label{Td*-eq}
\kappa(\Td(V))=\Td(V^\vee)=\Td(V)\cdot \ch(\det(V))^{-1}.
\end{equation}
In particular, the classes $\kappa(\Td(V))$ and $\kappa(\Td(V)^{-1})=\kappa(\Td(V))^{-1}$ 
belong to the rational lattice $J(H^*(X,\Q))$ 
(see Section \ref{rational-lattice-adm-sub-sec}).

Thus, the class $\kappa(\a\Td_X^{-1})$ is in $J(H^*(X,\Q))$. Hence, using \eqref{tau-formula} 
we get
$$\tau(y')=\kappa(\a\Td_X^{-1})\cdot\tau(y)=(-1)^p(2\pi i)^{q-p}\kappa(\a\Td_X^{-1})\cdot\ov{y}$$
Since the operator $\psi=\Pi\psi$ is compatible with the rational lattices (see Lemma \ref{quasihom-lattice-fun-lem}), we deduce 
$$\ov{\psi(y')}=\psi(\tau(y'))=(-1)^p(2\pi i)^{q-p}\psi(\kappa(\a\Td_X^{-1})\cdot\ov{y}).$$
Hence, by adjointness of the pair $(\psi,\phi)$,
\begin{equation}\label{main-y'-pairing-eq}
\lan\ov{\psi(y')},\psi(y')\ran_W=(-1)^p(2\pi i)^{q-p}\lan\kappa(\a\Td_X^{-1})\cdot\ov{y},\phi\psi(y')\ran_{D^b(X)}.
\end{equation}

Recall that $\phi\psi=\Phi_*\Psi_*=(\Phi\circ\Psi)_*$ is induced by the Fourier-Mukai functor 
with the kernel $K$, as in Lemma \ref{Koszul-comp-lem}.
Thus, by Lemma \ref{FM-HH-lem}, we have 
$$\phi\psi(y')=\tr_{12}(y'\ot \ch(K))=\int_{p_2} p_1^*(\Td_X\cdot y') \ch(K),$$ 
where 
$$\int_{p_2}:=(\int_X\ot\id):H^*(X\times X)\simeq H^*(X)\ot H^*(X)\to H^*(X).$$
Taking into account the relation \eqref{Td*-eq} we get
$$\phi\psi(y')=\int_{p_2} p_1^*(\ch(\om_X)^{-1}\kappa(\a)\cdot y) \ch(K).$$
Thus, we can rewrite \eqref{main-y'-pairing-eq} as follows:
\begin{align*}
&(-1)^p(2\pi i)^{p-q}\lan\ov{\psi(y')},\psi(y')\ran_W=
\lan \kappa(\a\Td_X^{-1})\cdot\ov{y},\phi\psi(y')\ran_{D^b(X)}=\int_X \a\cdot\kappa(\ov{y})\cdot\phi\psi(y')=\\
&(-1)^p\int_X \ov{y}\cdot\left(\a\cdot \int_{p_2}p_1^*(\ch(\om_X)^{-1}\kappa(\a)y)\cdot\ch(K)\right)=\\
&(-1)^p\int_{X\times X}p_2^*(\ov{y})\cdot p_1^*(y)\cdot \left(p_1^*(\ch(\om_X)^{-1}\kappa(\a))\cdot p_2^*\a\cdot \ch(K)\right).
\end{align*}
Note that the condition \eqref{main-homog-eq} simply means that
$$p_1^*(\ch(\om_X)^{-1}\kappa(\a))\cdot p_2^*\a\cdot \ch(K)\in H^{2D}(X\times X).$$
On the other hand, $p_2^*(\ov{y})\cdot p_1^*(y)\in H^{2j}(X\times X)$.
Thus, the above integral vanishes unless $2j+2D=4\dim X$. In other words, it is zero unless
$j=2\dim X-D$.
\ed

\subsection{Proof of Theorem \ref{dim-property-thm}}\label{dim-property-sec}

Let $W$ be a quasihomogeneous polynomial of degree $d$ as in Theorem \ref{dim-property-thm}, and let
$G$ be a finite group of diagonal symmetries of $W$, containing the exponential grading operator $J$.
Recall that for $\ga_1,\ldots,\ga_r\in G$ the maps 
$$\phi_g(\ga_1,\ldots,\ga_r):H(W_{\ga_1})^G\ot\ldots\ot H(W_{\ga_r})^G\to H^*(\SS,\C),$$
where $\SS=\SS_{g,\mu_d}(\ga_1,\ldots,\ga_r)$ is the moduli of $\Ga$-spin structures associated with $G$, 
giving the algebraic FJRW cohomological field theory, are obtained in the following way (see \cite[Sec.\ 5.1]{PV-CohFT}).
First, we consider the potential $W=W_{\ga_1}\oplus \ldots\oplus W_{\ga_r}$, and the map
$$\phi:HH_*(\MF_{\G_m}(W))\to H^*(\SS,\C)$$
defined as in Sec.\ \ref{purity-sec},
using a certain $\G_m$-equivariant Koszul matrix factorization $\{\a,\b\}$
of $-Z^*W$ on the total space of a vector bundle $p:\tot(A)\to\SS$, equipped with
a map $Z:\tot(A)\to\A^n$. 
The definition of the Koszul matrix factorization $\{\a,\b\}$ is rather involved (see \cite[Sec.\  4]{PV-CohFT}) and will not
be repeated here: for our purposes we only need to know that it is supported on the zero section in $\tot(A)$.
One difference from the framework of Sec.\ \ref{purity-sec} is that $\SS$ is not a variety, but a DM-stack. However,
there is still a natural map
\begin{equation}\label{Hoch-to-coh-map}
HH_*(\SS)\to H^*(\SS,\C)
\end{equation}
(see \cite[Eq.\ (5.6)]{PV-CohFT}), which we use to define $\phi$ with values in $H^*(\SS,\C)$.
Now the map $\phi_g$ is obtained by restricting $\Td(A)^{-1}\Td(B)\phi$ to the subspace
$$H(W_{\ga_1})^{\mu_d}\ot\ldots\ot H(W_{\ga_r})^{\mu_d}\sub H(W)^{\mu_d}.$$
(here we use the fact that $\Td(A)\Td(B)^{-1}=\Td(R\pi_*(\bigoplus_{j=1}^n\LL_j))$, where
$(\LL_\bullet)$ comes from a universal generalized spin-structure on $\SS$).

Actually, in \cite{PV-CohFT} we consider a bigger group $\Ga$  
and a $\Ga$-equivariant matrix factorization $\bP=(p,Z)_*\{\a,\b\}$ of $-W$ on $\A^n\times X$ to produce a map 
of $\C[G^*]$-modules
$$HH_*(\MF_{\Ga}(W))\to H^*(\SS,\C)\otimes \C[G^*],$$ 
where $G^*$ is the dual group to $G$. To get $\phi_g(\ga_1,\ldots,\ga_r)$ we specialize this map using the
evaluation at $1$ homomorphism $\C[G^*]\to \C$, compose the resulting map with a natural embedding
$$H(W_{\ga_1})^G\ot\ldots\ot H(W_{\ga_r})^G\to HH_*(\MF_{\Ga}(W))\ot_{\C[G^*]}\C,$$
and twist by $\Td(R\pi_*(\bigoplus_{j=1}^n\LL_j))^{-1}$.
It is easy to check that one gets the same map by passing to
$\G_m$-equivariant matrix factorizations and then applying the above procedure.

The map \eqref{Hoch-to-coh-map} is defined using 
a finite flat surjective morphism $\pi:X\to\SS$, where $X$ is a smooth projective variety (the existence of such maps
is a general fact about smooth proper DM-stacks over $\C$ with projective coarse moduli spaces---see 
\cite[Thm.\ 4.4]{Kr}, \cite[Thm.\ 2.1]{Kr-Vi}). In fact, \eqref{Hoch-to-coh-map} factors
through the pull-back map
$$\pi_*: HH_*(\SS)\to HH_*(X)\simeq H^*(X,\C),$$
followed by the degree preserving map $H^*(X,\C)\to H^*(\SS,\C)$.
Thus, it is enough to prove the required purity of dimension over $X$.
Taking into account the equality 
$$D_g(\ga_1,\ldots,\ga_r)=-\rk R\pi_*(\bigoplus_{j=1}^n\LL_j)=\rk B-\rk A,$$
we see that Theorem \ref{main-thm} would imply the dimension property \eqref{dimension-axiom} provided we
check that the polynomial $W_{\ga_1}\oplus\ldots\oplus W_{\ga_r}$ satisfies the assumptions of that Theorem,
i.e., each homogeneous polynomial $\wt{W_{\ga_i}}$ still has an isolated singularity at $0$. 
But this follows easily from Lemma \ref{isolated-lem} below.
\ed

\begin{lem}\label{isolated-lem} 
Let $W(x_1,\ldots,x_n)$ be a quasihomogeneous polynomial with an isolated singularity, where $\deg(x_i)=d_i>0$.
Let $I\sub [1,n]$ be the set of $i$ such that $d_i>1$. 
Then $\wt{W}(y_1,\ldots,y_n)=W(y_1^{d_1},\ldots,y_n^{d_n})$ still has an isolated singularity if and only
for every subset $J\sub I$ the restriction $W|_{\A^n_J}$ has an isolated singularity, 
where $\A^n_J\sub \A^n$ is the linear subspace given by $x_j=0$ for all $j\in J$.
\end{lem}

\Pf . We have $\pa_{y_i}\wt{W}(y)=d_iy_i^{d_i-1}\pa_{x_i} W(\varphi(y))$, where $\varphi(y_1,\ldots,y_n)=(y_1^{d_1},\ldots,y_n^{d_n})$. 
Thus, for $i\not\in I$ we have $\pa_{y_i}\wt{W}(y)=0$
if and only if $\pa_{x_i}W(\varphi(y))=0$. On the other hand,
for $i\in I$ we have $\pa_{y_i}\wt{W}(y)=0$ if and only if either $y_i=0$ or $\pa_{x_i}W(\varphi(y))=0$. This easily implies that
$y$ is a critical point of $\wt{W}$ if and only if $\varphi(x)$ is a critical point of $W|_{\A^n_J}$ for some subset $J\sub I$.
\ed


\section{APPENDIX: A compatibility involving the Grothendieck duality.}

Let $f:Y\to Z$ be a separated morphism of finite type between Noetherian schemes. 
We denote by $f^+:D^+(qcoh(Y))\to D^+(qcoh(Z))$ the extraordinary inverse image functor (see
\cite{Hart} where it is denoted by $f^!$).
Then for any $F\in D^b(Y)$ such that the support of $F$ is proper over $Z$, we have a canonical morphism
\begin{equation}\label{cF-eq}
c_{f,F}:F\to f^+f_*F.
\end{equation}
Indeed, this can be reduced to a similar map in the case when $f$ is proper: let $F=i_*F'$, where
$i:Y'\to Y$ is a closed subscheme, proper over $Z$, and $F'\in D^b(Y')$. Then we have $f_*F\simeq f'_*F'$, where $f'=f\circ i$. Since $f'$ is proper, we have a canonical map
$$c_{f',F'}:F'\to (f')^+f'_*F'\simeq i^+f^+f_*F.$$
By adjunction of $i_*$ and $i^+$ we get the required map 
$$F=i_*F'\to f^+f_*F.$$
Furthermore, for such $F$ and for any $G\in D^b(Z)$ the natural map
\begin{equation}\label{!-adj-map}
\Hom(f_*F,G)\to \Hom(f^+f_*F,f^+G)\to \Hom(F,f^+G),
\end{equation}
where the second arrow is induced by $c_F$, is an isomorphism (again this easily reduces to the case
when $f$ is proper).

The fact that $c_{f,F}$ does not depend on a choice of the subscheme $Y'$ follows from the compatibility
of the maps $c_{f,F}$ with compositions (for proper maps). Namely, for a morphism $g:Z\to T$ we have a commutative triangle
\begin{equation}\label{cF-compatibility}
\begin{diagram}
F&\rTo{c_{f,F}}& f^+f_*F\\
&\rdTo{c_{gf,F}}&\dTo_{f^+c_{g,f_*F}}\\
&&f^+g^+g_*f_*F
\end{diagram}
\end{equation}

The above picture extends to matrix factorizations. Namely, let us assume that we have the following situation:

\noindent
($\star$): {\it $Y$ and $Z$ are smooth $\G_m$-varieties admitting $\G_m$-equivariant ample line bundles; \break $f:Y\to Z$ is a smooth $\G_m$-equivariant
morphism; $W$ is a function on $Z$ of weight $d>0$ with respect to the $\G_m$-action, which is not a zero divisor.}

\noindent
Then canonical morphism \eqref{cF-eq} can be constructed for $F\in\MF_{\G_m}(Y,f^*W)$, with proper support over $Z$, using
Corollary \ref{duality-cor}, and the compatibility \eqref{cF-compatibility} still holds provided $g$ is also smooth 
$\G_m$-equivariant and $W=g^*W'$. 

Next, let us consider the fibered product $Y\times_Z Y$ with its two projections $p_1,p_2:Y\times_Z Y\to Y$, and let
$\de:Y\to Y\times_Z Y$ be the diagonal embedding. Let us also set
$$\pi=f\circ p_1=f\circ p_2: Y\times_Z Y\to Z.$$
Then for any $F\in D^b(Y)$ we have a canonical morphism
\begin{equation}\label{aF-map}
\a_F:\de_*F\to p_2^+F
\end{equation}
on $Y\times_Z Y$, which corresponds by adjunction to the identity map
$$F\to \de^+p_2^+F\simeq F.$$
Equivalently, it corresponds by adjunction to the identity map 
$$F\simeq p_{2*}\de_*F\to F$$
(note that $\de_*F$ is supported on the diagonal which is proper over $Y$).
Exchanging the roles of the factors in $Y\times_Z Y$ we get canonical morphisms
$$\a'_F:\de_*F\to p_1^+F.$$

As before, we can define similar morphisms for 
$F\in\MF_{\G_m}(Y,f^*W)$ assuming the situation ($\star$).

We will need the following properties of the maps $\a$ and $\a'$.
We set $\DD_f=f^+\OO_Z\simeq \om_f[\dim Y-\dim Z]$.

\begin{lem}\label{dual-base-change-lem} 
(i) In the above situation, assuming that $F\in D^b(Y)$ has the support that is proper over $Z$, 
we have a commutative triangle
\begin{diagram}
F\simeq p_{1*}\de_*F&\rTo{p_{1*}(\a_F)}&p_{1*}(p_2^+F)\\
&\rdTo{c_{f,F}}&\dTo{\th}\\
&&f^+f_*F
\end{diagram}
where 
$\th$ is the base change map, which corresponds via the isomorphism \eqref{!-adj-map} for $p_1$ to the map
$$p_2^+F\rTo{p_2^+c_{f,F}} p_2^+f^+f_*F\simeq p_1^+f^+f_*F$$
(note that the map from the support of $p_2^+F$ to $Y$, induced by $p_1$, is proper).
The similar assertion holds for $F\in\MF_{\G_m}(Y,f^*W)$, with proper support over $Z$, assuming the situation $(\star)$.

\noindent
(ii) For $F, G \in D^b(Y)$, or, assuming the situation $(\star)$, for $F\in\Per_{\G_m}(Y)$, $G\in\MF_{\G_m}(Y,f^*W)$,  the diagram
\begin{diagram}
\de_*(F)\ot p_2^*G&\rTo{\a_F\ot\id}&p_2^+F\ot p_2^*G\\
\dTo{\sim}&&\dTo{\sim}\\
\de_*(F\ot G)&\rTo{\a_{F\ot G}}&p_2^+(F\ot G)
\end{diagram}
is commutative. The same property holds for the maps $\a'$.

\noindent
(iii) Under the natural identification $p_1^+\DD_f\simeq p_2^+\DD_f\simeq p_1^*\DD_f\ot p_2^*\DD_f$ one has
$$\a_{\DD_f}=\a'_{\DD_f}\in\Hom(\de_*\DD_f, p_1^*\DD_f\ot p_2^*\DD_f).$$
\end{lem}

\Pf . (i) Applying the compatibility \eqref{cF-compatibility} to the maps $f$, $p_1$ and the object $\de_*F$ we get that
$c_{\pi,\de_*F}=c_{fp_1,\de_*F}$ is equal to the composition
$$\de_*F\to p_1^+p_{1*}\de_*F=p_1^+F\rTo{p_1^+c_{f,F}}p_1^+f^+f_*F.$$
In other words, the map $c_{f,F}$ corresponds to $c_{\pi,\de_*F}$ under the adjunction isomorphism
\eqref{!-adj-map}. Thus, we have to show that $c_{\pi,\de_*F}=c_{fp_2,\de_*F}$ is equal to the composition
$$\de_*F\rTo{\a_F} p_2^+F\rTo{p_2^+c_{f,F}} p_2^+f^+f_*F\simeq p_1^+f^+f_*F.$$
But this immediately follows from the compatibility \eqref{cF-compatibility} applied to the maps $f$, $p_2$ and
the object $\de_*F$.

\noindent
(ii) First, let us consider the case of sheaves.
By reversing the direction of the isomorphism of the left vertical arrow and using the adjointness of $(\de_*,\de^+)$
we reformulate the required commutativity as showing that the following composition is the identity map:
\begin{equation}\label{de!-composition}
F\ot G\to \de^+\de_*(F\ot G)\rTo{\sim} \de^+(\de_*(F)\ot p_2^*G)\rTo{\a_F}\de^+(p_2^+F\ot p_2^*G)\rTo{\sim}
\de^+p_2^+(F\ot G)\simeq F\ot G,
\end{equation}
where the first arrow is the adjunction map.
Now we use the following standard compatibility of the canonical morphisms $t_{f,A,B}:f^+(A)\ot f^*(B)\to f^+(A\ot B)$
with the projection formula: for a proper map $f$ the square
\begin{diagram}
f^+f_*(A\ot f^*B)&\rTo{\sim}&f^+(f_*A\ot B)\\
\uTo{c_{f,A\ot f^*B}}&&\uTo{t_{f,f_*A,B}}\\
A\ot f^*B&\rTo{c_{f,A}\ot\id}& f^+f_*A\ot f^*B
\end{diagram}
with the top horizontal arrow induced by the projection formula, is commutative.
Applying this to $f=\de$, $A=F$ and $g=p_2^*G$ we deduce the commutativity of the left square in the diagram
\begin{diagram}
\de^+\de_*(F\ot G)&\rTo{\sim}&\de^+(\de_*F\ot p_2^*G)&\rTo{}&\de^+(p_2^+F\ot p_2^*G)\\
\uTo{}&&\uTo{}&&\uTo{}\\
F\ot G&\rTo{}& \de^+\de_*F\ot \de^*p_2^*G&\rTo{}&\de^+p_2^+F\ot \de^*p_2^*G
\end{diagram}
Note that the second square is commutative by the functoriality of the map $\de^+A\ot\de^*B\to\de^+(A\ot B)$
in $A$. It follows that the composition of the first three arrows in \eqref{de!-composition} is equal to the map
\begin{equation}\label{de!-composition2}
F\ot G\simeq \de^+p_2^+F\ot \de^*p_2^*G\to \de^+(p_2^+F\ot p_2^*G).
\end{equation}

Next, for composable arrows $f$ and $g$ we have a commutative diagram
\begin{diagram}
(fg)^+(A)\ot (fg)^*(B)&\rTo{t_{fg,A,B}}&(fg)^+(A\ot B)\\
\dTo{\sim}&&\dTo{\sim}\\
g^+f^+(A)\ot g^*f^*(B)&\to g^+(f^+(A)\ot f^*(B))\to& g^+f^+(A\ot B)
\end{diagram}
Applying this for $f=p_2$, $g=\de$, $A=F$ and $B=G$ we deduce that the composition of
\eqref{de!-composition2} with the last arrow in \eqref{de!-composition},
$$\de^+(p_2^+F\ot p_2^*G)\rTo{\sim}\de^+p_2^+(F\ot G)\simeq F\ot G,$$
is the identity map of $F\ot G$.

The case of matrix factorization reduces to the case of ($\G_m$-equivariant)
sheaves using the equivalences with the ($\G_m$-equivariant) singularity categories.
Namely, let $Z_0\sub Z$ be the hypersurface of zeros of $W$, $Y_0=f^{-1}(Z_0)$.
Note that the hypersurface of zeros of $\pi^*W$ is
$$\pi^{-1}(Z_0)=Y_0\times_{Z_0} Y_0.$$
Now $G$ corresponds to an object of $D_{\Sg,\G_m}(Y_0)$, while
the commutative diagram lives in the category of matrix factorizations of $\pi^*W$, which is equivalent
to the category $D_{\Sg,\G_m}(Y_0\times_{Z_0} Y_0)$.
Now we observe that the functors $p_2^*$ and $p_2^+$ from matrix factorizations of $f^*W$ to those
of $\pi^*W$ correspond to the similar functors
$$p_2^*, p_2^+: D_{\Sg,\G_m}(Y_0,f^*W)\to D_{\Sg,\G_m}(Y_0\times_{Z_0} Y_0).$$
The operation of tensoring a matrix factorization with an object $P$ of the perfect derived category
corresponds for the category of singularity to the operation of tensoring with the restriction of $P$
to the zero locus of the potential (we apply this for $f^*W$ and for $\pi^*W$).
Finally, we use the fact that 
$$\de_*F|_{\pi^{_1}(Z_0)}\simeq \ov{\de}_*(F|_{Y_0}),$$
where $\ov{\de}:Y_0\to Y_0\times_{Z_0} Y_0$ is the diagonal,
that follows from the base change formula.

\noindent
(iii) Using the definition this reduces to checking the equality of the maps
$$\OO_Y \rTo{\sim} \de^+p_1^+f^+\OO_Z \ \ \text{ and}$$
$$\OO_Y \rTo{\sim} \de^+p_2^+f^+\OO_Z$$
under the identification $\pi^+\OO_Z\simeq p_1^+f^+\OO_Z\simeq p_2^+f^+\OO_Z$.
This reduces to the commutativity of the diagram with standard isomorphisms
\begin{diagram}
(fgh)^+&\rTo{\sim}&h^+(fg)^+\\
\dTo{\sim}&&\dTo{\sim}\\
(gh)^+f^+&\rTo{\sim}&h^+g^+f^+
\end{diagram} 
applied to the triples of morphisms $(f,p_1,\de)$ and $(f,p_2,\de)$.
\ed

Next, we assume that our map $f$ fits into a diagram
\begin{diagram}
&&Y\\
&\ldTo{f}&&\rdTo{p}\\
Z&&&& X
\end{diagram}
and that we are given $P\in D^b(Y)$, with proper support. Let us set $Q=P^\vee\otimes \DD_f$.
Note that we have a natural isomorphism $p_2^+\OO_Y\simeq p_1^*\DD_f$.
Hence, we obtain a natural map on $Y\times_Z Y$,
\begin{equation}\label{map-of-kernels2}
\begin{array}{l}
\wt{\varphi}:\de_*\OO_Y\to \de_*(P^\vee\ot P)\simeq \de_*\OO_Y\ot p_1^*P^\vee\ot p_2^*P\rTo{\a_{\OO}\ot \id} 
p_1^*\DD_f\ot p_1^*P^\vee\ot p_2^*P\simeq \\
p_1^*(\DD_f\ot P^\vee)\ot p_2^*P=p_1^*Q\ot p_2^*P,
\end{array}
\end{equation}
where $\de:Y\to Y\times_Z Y$ is the relative diagonal and $\a_\OO$ is the map \eqref{aF-map}.
Let $p_{XX}:Y\times_Z Y\to X\times X$ be the map with the
components $(pp_1,pp_2)$. Applying $p_{XX,*}$ to the above map and using the natural map $\OO_X\to p_*\OO_Y$
we get a canonical morphism
\begin{equation}\label{map-of-kernels}
\varphi:\De_*\OO_X\to \De_*p_*\OO_Y\simeq p_{XX,*}\de_*\OO_Y\to p_{XX,*}(p_1^*Q\ot p_2^*P),
\end{equation}
where $\De=\De_X:X\to X\times X$ is the diagonal map. 

We also consider an analogous construction for matrix factorizations in the situation ($\star$), where
$P$ is an object of $\MF_{\G_m}(Y,-f^*W)$ with proper support. We then view $Q$ as an object of $\MF_{\G_m}(Y,f^*W)$,
and the analog of the map \eqref{map-of-kernels2} can be constructed in
$\MF_{\G_m,d}(Y\times_Z Y,0)$. Assuming in addition that $X$ is smooth (so we can regard
$\De_*\OO_X$ as a perfect complex on $X\times X$), we get an analog of the map
$\varphi$ in $\MF_{\G_m,d}(X\times X,0)$.
Recall that in this situation we have functors
$$\wt{\Phi}_P:\MF_{\G_m}(Z,W)\to\MF_{\G_m,d}(X,0), \ \wt{\Psi}_Q:\MF_{\G_m,d}(X,0)\to \MF_{\G_m}(Z,W)$$
(see Sec.\ \ref{mf-ker-sec}). In the case of sheaves we also denote by $\wt{\Phi}_P$ and $\wt{\Psi}_Q$ the similar Fourier-Mukai  functors
between $D^b(Z)$ and $D^b(X)$. 

The main result of this Appendix is the following compatibility
(needed for the proof of Proposition \ref{mf-adjoint-prop}).

\begin{prop}\label{compatibility-prop} 
In the two situations described above (with sheaves and with matrix factorizations)
the map 
$$\Hom(\wt{\Psi}_{Q}(E),F)\to\Hom(E,\wt{\Phi}_P(F))$$
obtained via \eqref{adj-eq1} from the natural transformation
$\Id\to \wt{\Phi}_P\circ\wt{\Psi}_Q$, induced by \eqref{map-of-kernels},
is equal to the composition
\begin{equation}\label{main-adjunction-isom}
\begin{array}{l}
\Hom(f_*(p^*E\ot Q),F)\rTo{\sim}\Hom(p^*E\ot Q,f^*F\ot\DD_f)\simeq\Hom(p^*E,f^*F\ot P)\rTo{\sim}\\
\Hom(E,p_*(f^*F\ot P)).
\end{array}
\end{equation}
\end{prop}

\Pf . By definition, we have to prove that the composition \eqref{main-adjunction-isom} is equal to
the map
$$\Hom(\wt{\Psi}_{Q}(E),F)\rTo{\wt{\Phi}_P}\Hom(\wt{\Phi}_P\wt{\Psi}_Q(E),\wt{\Phi}_P(F))\to\Hom(E,\wt{\Phi}_P(F)),$$
where the second arrow is induced by \eqref{map-of-kernels}. Unraveling this leads to the following composition
\begin{align*}
&\Hom(f_*(p^*E\ot Q),F)\rTo{(1)} \Hom(p_{2*}p_1^*(p^*E\ot Q),f^*F)\rTo{(2)}\\
&\Hom(p^*E\ot p_{2*}(p_1^*Q\ot p_2^*P),f^*F\ot P)\rTo{(3)}
\Hom(E,p_*(f^*F\ot P)),
\end{align*}
with the intermediate maps given by
$$
(1):\Hom(f_*(p^*E\ot Q),F)\rTo{f^*}\Hom(f^*f_*(p^*E\ot Q),f^*F)\rTo{\sim} \Hom(p_{2*}p_1^*(p^*E\ot Q),f^*F),$$
where the second arrow is induced
by the base change isomorphism $f^*f_*\rTo{\sim} p_{2*}p_1^*$;
\begin{equation}\nonumber
\begin{array}{l}
(2):\Hom(p_{2*}p_1^*(p^*E\ot Q),f^*F)\rTo{\ot P}
\Hom(p_{2*}p_1^*(p^*E\ot Q)\ot P,f^*F\ot P)\simeq  \nonumber\\ 
\Hom(p_{2*}(p_1^*p^*E\ot p_1^*Q\ot p_2^*P),f^*F\ot P);
\end{array}
\end{equation}
and
\begin{align*}
&(3):\Hom(p_{2*}(p_1^*p^*E\ot p_1^*Q\ot p_2^*P),f^*F\ot P)\rTo{p_*}\\
&\Hom(p^X_{2*}p_{XX,*}(p_1^*p^*E\ot p_1^*Q\ot p_2^*P),p_*(f^*F\ot P))\simeq\\
&\Hom(p^X_{2*}\bigl(p^{X,*}_1E\ot p_{XX,*}(p_1^*Q\ot p_2^*P)\bigr),p_*(f^*F\ot P))
\rTo{\varphi}\\
&\Hom(p^X_{2*}(p^{X,*}_1E\ot \De_*\OO_X),p_*(f^*F\ot P))\simeq \Hom(E,p_*(f^*F\ot P)),
\end{align*}
where $p^X_i:X\times X\to X$, for $i=1,2$, are the projections. Here we used the natural identifications
$p_*p_{2*}\simeq p^X_{2*}p_{XX,*}$, $p_1^*p^*\simeq p_{XX}^*p^{X,*}_1$ and the projection formula for $p_{XX}$.

Let us set for brevity $R:=p_1^*Q\ot p_2^*P$ and $\wt{F}:=f^*F\ot P$. We claim that the map
$$(3'):\Hom(p_{2*}(p_1^*p^*E\ot R),\wt{F})\to \Hom(p^*E,\wt{F}),$$
that corresponds to (3) under the identification $\Hom(E,p_*(f^*F\ot P))\simeq \Hom(p^*E,f^*F\ot P)$,
is given simply by the composition 
$$\Hom(p_{2*}(p_1^*p^*E\ot R),\wt{F})\rTo{\wt{\varphi}}\Hom(p_{2*}(p_1^*p^*E\ot\de_*\OO_Y),\wt{F})\simeq
\Hom(p_{2*}\de_*(p^*E),\wt{F})\simeq\Hom(p^*E,\wt{F}),
$$
where the first arrow is induced by \eqref{map-of-kernels2}.
Indeed, first, using the definition of $\varphi$ we can rewrite (3) as the composition
\begin{align*}
&\Hom(p_{2*}(p_1^*p^*E\ot R),\wt{F})\rTo{p_*}\Hom(p^X_{2*}p_{XX,*}(p_1^*p^*E\ot R),p_*\wt{F})\rTo{\wt{\varphi}}\\
&\Hom(p^X_{2*}p_{XX,*}(p_1^*p^*E\ot \de_*\OO_Y),p_*\wt{F})\simeq
\Hom(p^X_{2*}(p^{X,*}_1E\ot p_{XX,*}\de_*\OO_Y),p_*\wt{F})\to\\
&\Hom(p^X_{2*}(p^{X,*}_1E\ot \De_*\OO_X),p_*\wt{F})\simeq\Hom(E,p_*\wt{F}).
\end{align*}
Now we observe that there is a commutative diagram
\begin{diagram}
\Hom(p_{2*}(p_1^*p^*E\ot R),\wt{F})&\rTo{p_*}&\Hom(p^X_{2*}p_{XX,*}(p_1^*p^*E\ot R),p_*\wt{F})\\
\dTo{\wt{\varphi}}&&\dTo{\wt{\varphi}}\\
\Hom(p_{2*}(p_1^*p^*E\ot \de_*\OO_Y),\wt{F})&\rTo{p_*}&\Hom(p^X_{2*}p_{XX,*}(p_1^*p^*E\ot \de_*\OO_Y),p_*\wt{F})\\
\dTo{\sim}&&\dTo{\sim}\\
\Hom(p_{2*}\de_*p^*E,\wt{F})&\rTo{p_*}&\Hom(p^X_{2*}p_{XX,*}\de_*p^*E,p_*\wt{F})\\
\dTo{\sim}&&\dTo{\sim}\\
\Hom(p^*E,\wt{F})&\rTo{p_*}&\Hom(p_*p^*E,p_*\wt{F})
\end{diagram}

Hence, our claim about the map $(3')$ follows from the commutativity of the diagram
\begin{diagram}
E&\rTo{\sim}&p^X_{2*}(p^{X,*}_1E\ot \De_*\OO_X)&\rTo{}& p^X_{2*}(p^{X,*}_1E\ot p_{XX,*}\de_*\OO_Y)\\
\dTo{}&&&&\dTo{\sim}\\
p_*p^*E&\rTo{\sim}&p^X_{2*}p_{XX,*}\de_*p^*E&\rTo{\sim}&p^X_{2*}p_{XX,*}(p_1^*p^*E\ot \de_*\OO_Y)
\end{diagram}
which is easy to check.

Our description of $(3')$ implies that we have a commutative diagram
\begin{equation}\label{b-ga-diagram}
\begin{diagram}
\Hom(p_{2*}(p_1^*(p^*E\ot Q)),f^*F)&\rTo{(3')\circ(2)}&\Hom(p^*E,f^*F\ot P)\\
\dTo_{\sim}^{\b}&&\dTo_{\sim}\\
\Hom(p_{2*}p_1^+(p^*E\ot Q),f^+F)&\rTo{\ga}&\Hom(p^*E\ot P^\vee\ot\DD_f,f^+F)
\end{diagram}
\end{equation}
Here $\b$ is the composition of the natural isomorphisms
\begin{align*}
&\b:\Hom(p_{2*}p_1^*(p^*E\ot Q),f^*F)\rTo{\sim}
\Hom(p_{2*}p_1^*(p^*E\ot Q)\ot \DD_f,f^*F\ot\DD_f)\rTo{\sim}\\
&\Hom(p_{2*}p_1^+(p^*E\ot Q),f^+F),
\end{align*}
where in the second isomorphism we use the identification $p_1^+\OO_Y\simeq p_2^*\DD_f$, and
$\ga$ is the composition
\begin{align*}
&\ga:\Hom(p_{2*}p_1^+(p^*E\ot Q),f^+F)\to\Hom(p_{2*}\bigl(p_1^+(p^*E\ot Q)\ot p_2^*P\ot p_2^*P^\vee\bigr),f^+F)\to\\
&\Hom(p^*E\ot P^\vee\ot\DD_f,f^+F),
\end{align*}
where the first arrow is induced by the evaluation map 
$ev_P:P\ot P^\vee\to \OO_Y$, while the second arrow is induced by the map
\begin{align}\label{kappa-id-map}
&p^*E\ot P^\vee\ot \DD_f\simeq p_{2*}\bigl(p_1^*p^*E\ot \de_*(P^\vee\ot\DD_f)\bigr)\rTo{p_{2*}(\id\ot\eps)}
p_{2*}(p_1^*p^*E\ot p_1^+Q\ot p_2^*P\ot p_2^*P^\vee)\simeq \nonumber\\
&p_{2*}\bigl(p_1^+(p^*E\ot Q)\ot p_2^*P\ot p_2^*P^\vee\bigr),
\end{align}
where $\eps$ is the composition
$$\eps:\de_*(P^\vee\ot\DD_f)\simeq
\de_*\OO_Y\ot p_2^*P^\vee\ot p_2^*\DD_f\rTo{\wt{\varphi}\ot\id} p_1^*Q\ot p_2^*P\ot p_2^*P^\vee\ot p_2^*\DD_f\simeq
p_1^+Q\ot p_2^*P\ot p_2^*P^\vee.$$

Next, let us consider the composed map
$$(1'):\Hom(f_*(p^*E\ot Q),F)\rTo{(1)}\Hom(p_{2*}p_1^*(p^*E\ot Q),f^*F)\rTo{\b}
\Hom(p_{2*}p_1^+(p^*E\ot Q),f^+F).
$$
Applying Lemma \ref{dual-base-change-lem}(i) to the object $p^*E\ot Q$ we get that the composition
$$p^*E\ot Q\simeq p_{2*}\de_*(p^*E\ot Q)\rTo{p_{2*}(\a_{p^*E\ot Q})}p_{2*}p_1^+(p^*E\ot Q)\to f^+f_*(p^*E\ot Q)$$
is just the map \eqref{cF-eq}. 
This implies the commutativity of the following triangle 
\begin{equation}\label{ga'-diagram}
\begin{diagram}
\Hom(f_*(p^*E\ot Q),F)&\rTo{(1')}&\Hom(p_{2*}p_1^+(p^*E\ot Q),f^+F)\\
&\rdTo{\sim}&\dTo{\ga'}\\
&&\Hom(p_*E\ot Q,f^+F)
\end{diagram}
\end{equation}
where the diagonal is the adjunction isomorphism and $\ga'$ is induced by
the map 
$$p^*E\ot Q\simeq p_{2*}\de_*(p^*E\ot Q)\rTo{p_{2*}(\a'_{p^*E\ot Q})} p_{2*}p_1^+(p^*E\ot Q).$$
Comparing the diagrams \eqref{b-ga-diagram} and \eqref{ga'-diagram} with the definition of the map
\eqref{main-adjunction-isom}, we see that our assertion would follow from the equality 
$\ga'=\ga$, which in turn would be implied by the commutativity of the diagram
\begin{diagram}
p^*E\ot P^\vee\ot\DD_f&\rTo{}&p_{2*}\bigl(p_1^+(p^*E\ot Q)\ot p_2^*P\ot p_2^*P^\vee\bigr)\\
\dTo{=}&&\dTo{ev_{p_2^*P}}\\
p^*E\ot Q&\rTo{p_{2*}(\a'_{p^*E\ot Q})}&p_{2*}p_1^+(p^*E\ot Q)
\end{diagram}
where 
the top horizontal arrow is \eqref{kappa-id-map}. Note that this diagram is obtained by applying the functor $p_{2*}$ to the diagram
\begin{diagram}
p_1^*p^*E\ot \de_*(P^\vee\ot\DD_f)&\rTo{\id\ot\eps}&p_1^+(p^*E\ot Q)\ot p_2^*P\ot p_2^*P^\vee\\
\dTo{\sim}&&\dTo{ev_{p_2^*P}}\\
\de_*(p^*E\ot Q)&\rTo{\a'_{p^*E\ot Q}}&p_1^+(p^*E\ot Q)
\end{diagram}
Unraveling the definition of $\eps$ and using Lemma \ref{dual-base-change-lem}(ii),
we see that the latter diagram is obtained by tensoring with $p_1^*p^*E$ from the diagram
\begin{equation}\label{last-diagram}
\begin{diagram}
\de_*\DD_f\ot p_2^*P^\vee&\rTo{}&p_2^+\DD_f\ot p_1^*P^\vee\ot p_2^*P\ot p_2^*P^\vee\\
\dTo{\sim}&&\dTo{ev_{p_2^*P}}\\
\de_*\DD_f\ot p_1^*P^\vee&\rTo{\a'_{\DD_f}\ot\id}p_1^+\DD_f\ot p_1^*P^\vee\simeq& p_2^+\DD_f\ot p_1^*P^\vee
\end{diagram}
\end{equation}
where the top arrow is the composition
$$\de_*\DD_f\ot p_2^*P^\vee\to
\de_*\DD_f\ot p_2^*P^\vee\ot p_2^*P\ot p_1^*P^\vee
\rTo{\a_{\DD_f}\ot\id} p_2^+\DD_f\ot p_2^*P^\vee\ot p_2^*P\ot p_1^*P^\vee.$$
Thus, it remains to prove the commutativity of \eqref{last-diagram}.
Using the commutative diagram
\begin{diagram}
\de_*\DD_f\ot p_2^*P^\vee\ot p_2^*P\ot p_1^*P^\vee&\rTo{\a_{\DD_f}\ot\id}&
p_2^+\DD_f\ot p_2^*P^\vee\ot p_2^*P\ot p_1^*P^\vee\\
\dTo{ev_{p_2^*P}}&&\dTo{ev_{p_2^*P}}\\
\de_*\DD_f\ot p_1^*P^\vee&\rTo{\a_{\DD_f}\ot\id}&p_2^+\DD_f\ot p_1^*P^\vee
\end{diagram}
we can rewrite the composition of the top arrow with the right vertical arrow in \eqref{last-diagram}
as the map
$$\de_*\DD_f\ot p_2^*P^\vee\to \de_*\DD_f\ot p_2^*P^\vee\ot p_2^*P\ot p_1^*P^\vee\rTo{ev_{p_2^*P}}\de_*\DD_f\ot p_1^*P^\vee\rTo{\a_{\DD_f}\ot\id} p_2^+\DD_f\ot p_1^*P^\vee.$$
Here the composition of the first two arrows coincides with the left vertical arrow in \eqref{last-diagram},
so we get the commutativity of the diagram like \eqref{last-diagram} but with $\a_{\DD_f}$ instead of $\a'_{\DD_f}$
in the bottom arrow. It remains to recall that $\a_{\DD_f}=\a'_{\DD_f}$ by Lemma \ref{dual-base-change-lem}(iii).
\ed

{\sc Department of Mathematics, University of Oregon, Eugene, OR 97405}

\end{document}